\documentclass[10pt,a4paper]{article}

\usepackage{amssymb,amsbsy,amsmath,amsfonts,amscd,mathrsfs}
\usepackage{graphicx}
\usepackage{color}
\usepackage[usenames,dvipsnames]{xcolor}
\usepackage{a4wide}
\usepackage{enumerate}

\newcommand{\proof}{\noindent {\bf Proof. }}  % COMMENTARE PER CMP
\newtheorem{remark}{Remark}  % COMMENTARE PER CMP
\newtheorem{lemma}{Lemma}       % COMMENTARE PER CMP

\newcommand{\llabel}[1]{{\label{#1}}}

\newcommand{\ffoot}[1]{}
\newcommand{\bi}{\begin{itemize}}
\newcommand{\ei}{\end{itemize}}
\newcommand{\bd}{\begin{description}}
\newcommand{\ed}{\end{description}}
\newcommand{\be}{\begin{enumerate}}
\newcommand{\ee}{\end{enumerate}}

\newcommand{\mm}{\boldsymbol{m}}
\newcommand{\beq}{\begin{equation}}
\newcommand{\eeq}{\end{equation}}
\newcommand{\bqn}{\begin{eqnarray}}
\newcommand{\eqn}{\end{eqnarray}}
\newcommand{\eqnn}{\nonumber\end{eqnarray}}
\newcommand{\eqnl}[1]{\llabel{#1}\end{eqnarray}}

\newcommand{\ba}[1]{\begin{array}{#1}}
\newcommand{\ea}{\end{array}}

\newcommand{\R}{\mathbb{R}}

\newcommand{\N}{\mathbb{N}}

\newcommand{\fine}{\end{document}}

\def \trait (#1) (#2) (#3){\vrule width #1pt height #2pt depth #3pt}
\def \qed{\hfill
        \trait (0.1) (6) (0)
        \trait (6) (0.1) (0)
        \kern-6pt   
        \trait (6) (6) (-5.9)
        \trait (0.1) (6) (0)
\medskip}

\newtheorem{Theorem}{\bf Theorem}[section]

\newtheorem{mcc}[Theorem]{\bf Corollary}
\newtheorem{Definition}[Theorem]{\bf Definition}
\newtheorem{mpr}[Theorem]{\bf Proposition}
\newtheorem{mproperty}[Theorem]{\bf Property}

\newcommand{\bt}{\begin{Theorem}}
\newcommand{\et}{\end{Theorem}}
\newcommand{\bp}{\begin{mpr}}
\newcommand{\ep}{\end{mpr}}
\newcommand{\bcc}{\begin{mcc}}
\newcommand{\ecc}{\end{mcc}}
\newcommand{\bex}{\begin{example}}
\newcommand{\eex}{\end{example}}

\newcommand{\bproperty}{\begin{mproperty}}
\newcommand{\eproperty}{\end{mproperty}}

\newcommand{\ec}{\end{mcc}}
\newcommand{\bdeff}{\begin{Definition}}
\newcommand{\edeff}{\end{Definition}}

\newcommand{\eproof}{\hfill $\square$}

\newcommand{\brem}{\begin{remark}}
\newcommand{\erem}{\end{remark}}

\newcommand{\lam}{\lambda}

\newcommand{\al}{\alpha}
\newcommand{\eps}{\varepsilon}

\newcommand{\om}{\omega}

\newcommand{\con}{{\cal C}}

\renewcommand{\H}{{\cal H}}

\newcommand{\mariospectrum}{separated discrete spectrum}
\newcommand{\mariomatrix}{conicity matrix}
\newcommand{\campopao}{non-mixing field}
\newcommand{\Xpao}{\mathcal{X}_P}
\newcommand{\Xpaot}{\mathcal{X}_{\mu}}

\newcommand{\uu}{{\bf u}}
\newcommand{\ww}{{\bf w}}
\newcommand{\vv}{{\bf v}}

\newcommand{\zz}{{\bf z}}

\newcommand{\Hdue}{\mathcal{W}}
\newcommand{\Heff}{H^{\varepsilon}_{\mathrm{eff}}}
\newcommand{\Ueff}{U^{\varepsilon}_{\mathrm{eff}}}

\newcommand{\AV}{\mathbf{A}}

\newcommand{\bth}{\boldsymbol{\Xi}}

\title{\LARGE {\bf
Approximate controllability via adiabatic techniques for the three-inputs controlled Schr\"odinger equation.}
}

\author{Francesca Carlotta Chittaro \thanks{ Aix Marseille Universit\'e, CNRS, ENSAM, LSIS UMR 7296, 13397 Marseille,  % FORMATO ARTICLE
France, and Universit\'e de Toulon, CNRS, LSIS UMR 7296, 83957 La Garde,  France {\tt francesca-carlotta.chittaro@univ-tln.fr}} 
\and Paolo Mason
 \thanks{
Laboratoire des Signaux et Syst\`emes (L2S, UMR 8506), CNRS - CentraleSupelec - Universit\'e Paris-Sud, 3, rue Joliot Curie, 91192, Gif-sur-Yvette, France {\tt mason@lss.supelec.fr}.} }

\date{}

\begin{document}
\maketitle

\begin{abstract}
We consider a system described by a controlled bilinear Schr\"odinger equation with three external inputs. We provide a constructive method to 
approximately steer the system from a given energy level to a superposition of energy levels corresponding to a given probability distribution.
The method is based on adiabatic techniques and works if the spectrum of the Hamiltonian
admits eigenvalue intersections, with respect to variations of the controls, and if the latter are conical. We provide  sharp estimates of the relation
between the error and the controllability time, and we show how to improve these estimates by selecting special control paths.
\end{abstract}

\section{Introduction}

A typical issue in quantum control concerns the controllability of the bilinear Schr\"odinger equation
\bqn
i\frac{d\psi}{dt}=\left(H_0+\sum_{k=1}^m u_k(t)H_k\right)\psi(t),
\label{hamiltonianam}
\eqn
where $\psi$ belongs to the Hilbert sphere  %$\mathbf{S}$ 
of a { (finite or infinite dimensional)} complex separable Hilbert space ${\cal H}$ and $H_0,\ldots, H_m$  
are self-adjoint operators on ${\cal H}$.
Here $H_1,\ldots, H_m$ represent  the action of external fields on the system, whose strength is given by
the scalar-valued controls $u_1,\ldots,u_m$, while $H_0$ describes the uncontrolled dynamics of the system.

The controllability problem aims at establishing whether, for every pair of states $\psi_0$ and $\psi_1$ in %$\mathbf{S}$
the Hilbert sphere, there exist  controls $u_k(\cdot)$ and a time $T$ such that the solution of~\eqref{hamiltonianam} with initial condition $\psi(0)=\psi_0$ satisfies $\psi(T)=\psi_1$. 

While the case where $\mathcal{H}$ is a finite dimensional Hilbert space has been widely understood~\cite{albertini,dalessandro-book}, in the infinite dimensional case the answer is far from being given. In particular,
negative results have been proved when ${\cal H}$ is infinite-dimensional (see~\cite{bms,turinici}).
Hence one has to look for weaker controllability properties as, for instance, approximate controllability (see for instance~\cite{schNOI2,schNOI1,mirrahimi-continuous,nersesyan}), or controllability  between subfamilies of states  (in particular
the eigenstates of $H_0$, which are the most relevant physical states) or in more regular subspaces of square-integrable functions (see~\cite{beauchard-coron,camillo}).

While the above mentioned works are essentially obtained by means of non-constructive arguments, the purpose of this paper is to propose a method that permits to explicitly 
select control inputs steering the system from the initial state to an arbitrarily small neighborhood of the given target state.
Adiabatic theory and conical intersections between eigenvalues constitute the main tools of the control strategy we propose in this paper.

Roughly speaking, the adiabatic theorem (see~\cite{bofo,nenciu,teufel}) states 
that the occupation probabilities associated with the energy levels of a time-dependent Hamiltonian $H(\cdot)$ are approximately preserved along the evolution given by $i\dot\psi(t)=H(t)\psi(t)$, provided that $H(\cdot)$ varies slowly enough. This result works whenever the energy levels (i.e. the eigenvalues of $H(\cdot)$) are isolated for every $t$. On the other hand, if two eigenvalues intersect, and
provided that $H(\cdot)$ is smooth enough, the passage through the intersections may determine (approximate) exchanges of the corresponding occupation 
probabilities (see~\cite[Corollary 2.5]{teufel}  and \cite{hagedorn}). For these reasons, adiabatic methods are largely used in quantum control to induce population transfers (see for instance the techniques known as Stimulated Raman Adiabatic Passage
(STIRAP), Stark-chirped rapid adiabatic passage (SCRAP)) and to prepare superposition states \cite{cogwheel}.
The applications of adiabatic methods in quantum control, as a tool for obtaining controllability results, have already been exploited in previous papers 
(see for instance~\cite{Boscain_Adami,adia-nostro,rouchon_a,jauslin}).
The general idea is to use slowly varying controls, taking advantage of the adiabatic theorem, and ``climb'' the energy levels through the conical intersections.

Related to the present paper is a method %that permits to attain (approximately) a prescribed distribution of probability starting from an eigenstate has been 
recently developed in~\cite{adia-nostro} in the case $m=2$ and for self-adjoint Hamiltonians with real matrix elements. 
It exploits a generalization of \cite[Corollary 2.5]{teufel} stating that it is possible to arbitrarily recombine the probability weights associated with two subsequent energy levels by following (slowly) a suitable control path passing through a conical intersection between them. The control strategy of~\cite{adia-nostro} %, that will be pursued in this paper under different assumptions, is applicable 
applies whenever a part of the spectrum of the Hamiltonian operator is uniformly separated from the rest of the spectrum (as a function of the control parameter), is discrete and each pair of subsequent eigenvalues intersect in a conical intersection.
When there exists such a portion of the spectrum, called \textit{separated discrete spectrum}, %and under additional technical assumptions on the Hamiltonian, 
this control strategy permits to attain (approximately) a state having  a prescribed distribution of probability 
(relative to the energy levels of the separated discrete spectrum) starting from an eigenstate. %/given energy level. 
In particular this entails a controllability property, that we call \textit{spread controllability}, which, although weaker than the usual approximate controllability property, is more practical. Note indeed that the relative phases between pairs of components in the eigenbasis decomposition are essentially uncontrollable since 
they evolve according to the gaps between the corresponding energy levels.
%is quite important in an experimental setting.....}
Furthermore, notice that this method allows us to control the population inside some portion of the discrete spectrum, 
if well separated from the rest, even in the presence of continuous spectrum, unlike many other classical methods.

Concerning the precision of the method, an application of the adiabatic theorem together with~\cite[Corollary 2.5]{teufel} shows that the maximal error is of the order of the square root of the control speed. On the other hand in~\cite{adia-nostro} it was shown that   the precision of the transfer may be remarkably improved if one follows some special paths in the space of controls; namely, such paths permit to attain a state with a prescribed probability distribution with an error of the order of the control speed. From a practical point of view this means that, to guarantee a given precision, one may significantly reduce the duration of the process, whose extent constitutes one of the main disadvantages of the implementation of adiabatic techniques.

The purpose of this paper is to adapt the control strategy introduced in~\cite{adia-nostro} to the general case of self-adjoint Hamiltonians, assuming that three controlled Hamiltonians are employed, and to select the control paths that allow to improve the precision of the process as explained above. Preliminary results in this sense were discussed in~\cite{cdc-maui}. Notice that the chosen setting is quite natural, since it is well-known, for Hermitian matrices or within spaces of self-adjoint operators satisfying particular transversality conditions, that the set of operators admitting multiple eigenvalues is a submanifold  of codimension three (see e.g.~\cite{agrachev,teytel,von_neumann-wigner}). Moreover  conical intersections do not constitute a pathological phenomenon since, as shown in  Appendix~\ref{gen}, all eigenvalue intersections are generically conical in the finite dimensional case and in some physically relevant infinite dimensional models. Conical intersections are  also structurally stable 
with respect to variations of the Hamiltonian operator, as shown in Theorem~\ref{struct}. Concerning the relation between conical intersections and controllability properties of the bilinear Schr\"odinger equation, let us finally mention the main results of the recent paper~\cite{ugauthier}: if all subsequent energy levels of the Hamiltonian are connected by means of conical intersections then the system is approximately controllable and, in the finite dimensional case, it is even exactly controllable.
Notice that these results have not been obtained by adiabatic techniques, although,  as shown in Section~\ref{fin rem}, it is not difficult to recover approximate controllability of the system by extending the results of this paper (in a non-constructive way).
 
The structure of the paper is the following. In Section~\ref{preliminary} we introduce the notations used throughout the paper, the main assumptions and definitions, and we adapt the classical  statement of the adiabatic theorem to our setting. In Section~\ref{s-conic} we discuss some properties of conical intersections and related results that allow to propose the basic control strategy. In Section~\ref{una-a-caso} we define special paths and, by means of a series of technical results, we show that they can be included in the control algorithm in order to improve its performance. As a byproduct, we get a structural stability result concerning conical intersections. In Section~\ref{fin rem} we briefly mention some extensions of the control strategy and of the controllability results obtained earlier.  Appendix~\ref{reg proj} gathers the technical results concerning the regularity of the spectrum and of the spectral projections that are needed throughout the paper, while Appendix~\ref{gen} 
discusses the 
genericity of conical intersections in the finite and infinite dimensional cases.

%%%%%%%%%%%%%%%

\section{Notations and preliminary results}
\label{preliminary}

We start this section by introducing the notations that will be used in the rest of the paper.\\
For a function $f(\cdot)$ of a real parameter $s$, we use the following notation for its right and left limits at $s_0$:
\[
f(s_0^{\pm})=\lim_{s\to s_0^{\pm}} f(s). 
\]
Moreover we say that $f(s)=o(s^k)$ if $\lim_{s\to 0}\frac{f(s)}{s^k}=0$.\\
Whenever $\gamma(s),\ s\in [s_1,s_2]$ is a curve on $\R^3$ and $Q(\cdot)$ is a function of $\vv\in\R^3$ then, with abuse of notations, we denote 
by $\dot Q(\gamma(r))$ the derivative of the composition $Q(\gamma(\cdot))$ computed at $r$, that is 
$\dot Q(\gamma(r)):=\frac{d}{ds}Q(\gamma(s))|_{s=r}=\frac{dQ}{d\vv}(\gamma(r))\frac{d\gamma}{ds}(r)$.
Similarly,  $Q^{(l)}(\gamma(r)):=\frac{d^l}{ds^l}Q(\gamma(s))|_{s=r}$.\\ 
The scalar product of two elements $\psi_1,\psi_2$  in the Hilbert state space 
is denoted by $\langle \psi_1,\psi_2\rangle$, while the scalar product of two vectors $\ww_1,\ww_2$ in any other euclidean space 
is denoted by $\ww_1\cdot \ww_2$. 
Analogously, the norm in the two cases is denoted respectively by $\|\cdot\|$ and $|\cdot|$.\\
For a given vector $\vv$ or matrix $A$ the respective transpose is denoted by $\vv^T$ and $A^T$. The inverse of the transpose of an invertible square matrix $A$ is denoted with $A^{-T}$.\\
Given a vector $\vv=(v_1,v_2,v_3)\in \mathbb{C}^3$, 
we denote its complex conjugate $(v_1^*,v_2^*,v_3^*)$ by $\vv^*$
and its real and imaginary parts respectively by
\[
\mathfrak{Re}(\vv)=(\mathfrak{Re}(v_1),\mathfrak{Re}(v_2),\mathfrak{Re}(v_3))
\qquad
\mathfrak{Im}(\vv)=(\mathfrak{Im}(v_1),\mathfrak{Im}(v_2),\mathfrak{Im}(v_3)).
\]
The symbol $\mathrm{id}$ is used to denote the identity operator on a vector space which is specified at each occurrence, whenever not clear from the context.\\

\subsection{General setting}
\label{definizionieris}

Let $\mathcal{H}$ be a separable complex  Hilbert space with norm $\|\cdot\|$;
let us introduce the following notion of relative boundedness between operators:

\bdeff[$A$-smallness and $A$-boundedness] Let $A,B$ two densely defined operators with domains $\mathcal{D}(A)\subset\mathcal{D}(B)$. We say that $B$ is $A$-bounded  
if there exist $a,b>0$ such that
$\|B\psi\|\leq a \|A\psi\|+b\|\psi\|$
for every $\psi\in\mathcal{D}(A)$.
$B$ is said to be $A$-small 
if for every $\alpha>0$ there exists $\beta>0$ such that  $\|B\psi\|\leq \alpha\|A\psi\|+\beta\|\psi\|$ for  every $\psi\in\mathcal{D}(A)$. 
(The latter notion is called \emph{infinitesimal smallness with respect to $A$} in~\cite{reed_simon_2}.)\label{ksb}
\edeff
Given a self-adjoint operator $A$ on  $\mathcal{H}$, for every $A$-bounded operator $B$ 
we define its norm with respect to $A$ as 
\begin{equation} \label{norm}
\|B\|_A= \sup_{\psi\in\mathcal{D}(A)} \frac{\|B\psi\|}{\|A\psi\|+\|\psi\|}.\end{equation}
This provides a norm in the space
 $\mathcal{L}(\mathcal{D}(A),\mathcal{H})$ of continuous linear operators from $\mathcal{D}(A)$ (endowed with the graph norm of $A$) to $\mathcal{H}$.

\medskip

We consider the Hamiltonian 
\[
H(\uu) = H_0 + u_1 H_1+u_2 H_2+ 
u_3 H_3,
\]
with $\uu = (u_1, u_2,u_3) \in \R^3$, and where $H_i,\ i=0,\ldots,3$ satisfy the following assumption:

\smallskip
\noindent
{\bf (H0)} $H_0$ is a self-adjoint operator on a separable complex Hilbert space ${\cal H}$, and $H_i$ are $H_0$-small self-adjoint operators on ${\cal H}$ for $i=1,2,3$.

%\brem
Under assumption {\bf (H0)}, \cite[Theorem X.12]{reed_simon_2} guarantees that $H(\uu)$ is self-adjoint with domain $\mathcal{D}(H_0)$.
Moreover, it is easy to see that for every $\uu$, $H_0$ is $H(\uu)$-bounded, and therefore $H_i$ is $H(\uu)$-small, for every $i=1,2,3$,	 
with constants $a,b$ (as in Definition~\ref{ksb}) that depend continuously on $\uu$. 
%\erem

Schr\"odinger Hamiltonians are typical Hamiltonian operators describing quantum phenomena and can be represented 
in the form $-\Delta+V$ on the Hilbert space $L^2(\Omega)$, where $\Omega$ is a domain of $\R^n$,  $\Delta$ is the Laplacian on $\Omega$ (with Dirichlet or Neumann boundary conditions) and $V:\Omega\to\R$ has to be interpreted as a multiplicative operator on $L^2(\Omega)$. In particular such Hamiltonian operators are unbounded operators.
In this context Hypothesis {\bf (H0)} is thus intended to describe an Hamiltonian operator of the previous form that can be controlled by means of three external inputs so that $H_0=-\Delta+V_0$ and $H_i=V_i$ for some multiplicative operators  $V_i$, for $0\leq  i\leq 3$.
 
Finite dimensional representations of quantum systems are also common, for instance in the description of spin systems. In this case the Hamiltonian operator $H(\uu)$ is a Hermitian matrix.  Consider for instance the case of a spin-$1/2$ particle immersed in a controlled magnetic field. In this case, $H_i$ are the Pauli matrices, and
the controls are the components of the magnetic field.

\smallskip

The dynamics of the quantum systems we consider are described by the time-dependent Schr\"odinger equation
\bqn
i\frac{d\psi}{dt}=H(\uu(t))\psi(t).
\label{hamiltonianabold}
\eqn
Such an equation has mild solutions under hypothesis {\bf (H0)}, $\uu(\cdot)$ piecewise $\con^1$ and with an initial condition in the domain of $H_0$ 
(see e.g.~\cite[Theorem X.70]{reed_simon_2} and~\cite{bms}).

\medskip

We are interested in controlling~\eqref{hamiltonianabold} inside some portion of the discrete spectrum of $H(\uu)$. 
Since we use adiabatic techniques, 
some spectral gap condition is needed: 

\smallskip
\noindent
{\bf (H1)} There exist a domain $\omega$ in $\R^3$, a map $\Sigma$ defined on $\omega$ that associates with each $\uu\in\omega$ a subset  $\Sigma(\uu)$ of the discrete spectrum of $H(\uu)$, and two continuous functions $f_1,f_2:\omega\to\R$ such that
\bi
\item $f_1(\uu) < f_2(\uu)$ and $\Sigma(\uu)\subset [f_1(\uu) , f_2(\uu)]\qquad\forall \uu\in\omega$.
\item there exists $\Gamma>0$ such that
\[\inf_{\uu\in\omega} \inf_{\lambda \in \mathrm{Spec}(H(\uu))\setminus \Sigma(\uu)} \mathrm{dist}(\lambda,[f_1(\uu), f_2(\uu)])  )>\Gamma.
\]
\ei
In this case we say that $\Sigma(\uu)$ is a \mariospectrum.

\medskip
\noindent{\bf Notation\ } From now on we label  the eigenvalues  belonging to a \mariospectrum\ $\Sigma(\uu)$ in such a way that %we can write 
$\Sigma(\uu)=\{\lambda_0(\uu) ,\ldots,\lambda_k(\uu)\}$, where  
$\lambda_0(\uu)\leq \cdots \leq\lambda_k(\uu)$ are counted according to their multiplicity (note that the separation of $\Sigma$ from the rest of the spectrum guarantees that $k$ is constant).
Moreover we denote by  $\phi_0(\uu) , \ldots , \phi_k(\uu)$ an orthonormal family of eigenstates corresponding to $\lambda_0(\uu), \ldots , \lambda_k(\uu)$.  Notice that in this notation $\lambda_0$ does not need to be the ground state of the system.

\smallskip

\begin{figure}
\begin{center}
\includegraphics[width=0.85\textwidth]{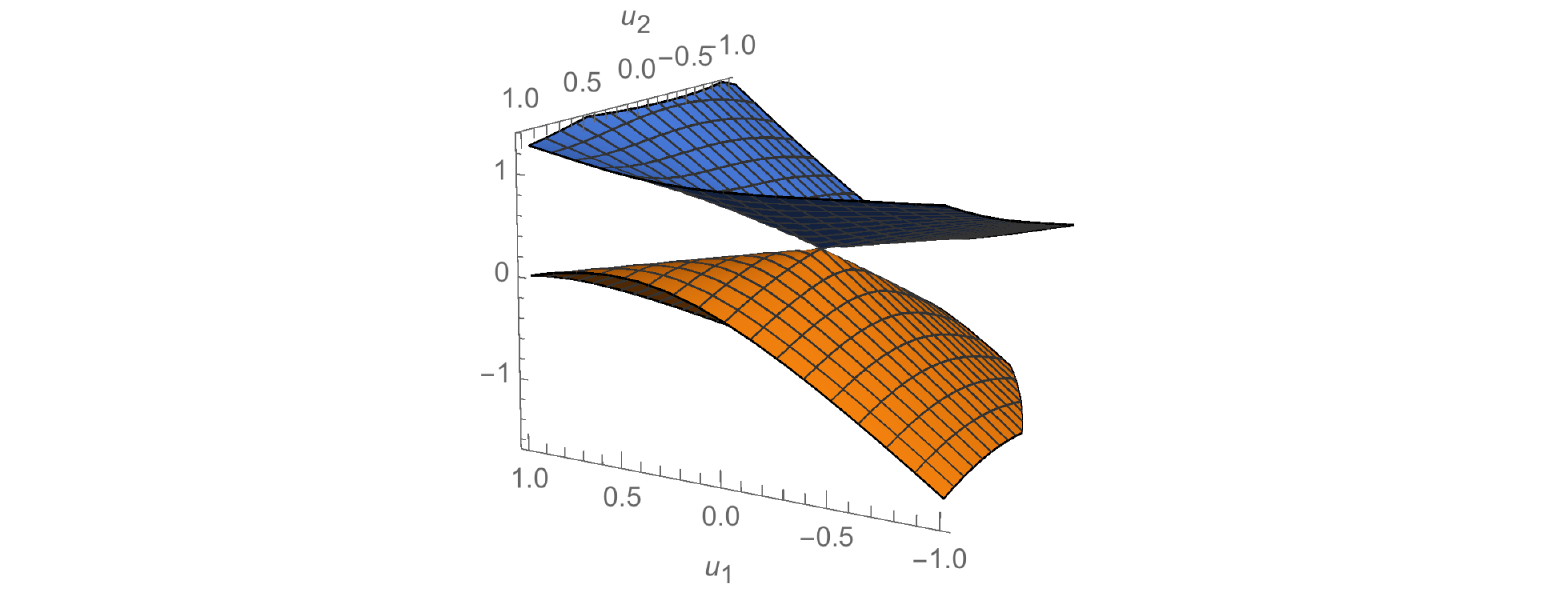}
\caption{A conical intersection for a controlled Hamiltonian with two inputs.}\label{f-conical}
\end{center}
\end{figure}

Our techniques rely on the existence of conical intersections between the eigenvalues, which constitute a well studied phenomenon in molecular physics 
(see for instance~\cite{bofo,lasser-fermanian,hagedorn,lasser,diabolical}). 
In this paper we will adopt the following definition, consistent with the one already given in~\cite{adia-nostro} for the two-inputs case 
(Figure~\ref{f-conical} shows a conical intersection in this latter setting).

\bdeff\label{conical}
Let $H(\cdot)$ satisfy hypothesis {\bf (H0)}. We say that $\bar\uu\in\R^3$ is a \emph{conical intersection} between two subsequent eigenvalues $\lam_{j}$ and $\lambda_{j+1}$ if $\lam_{j}(\bar \uu) = \lam_{j+1}(\bar \uu)$ has multiplicity two and there exists a constant $c>0$ such that for any unit vector $\mathbf{v}\in \R^3$ and $t>0$ small enough we have that 
\begin{equation} \label{formcono}
\lam_{j+1}(\bar \uu+t\mathbf{v})-\lam_{j}(\bar \uu+t\mathbf{v}) > ct\,.\end{equation}
\edeff
A discussion on the occurrence of conical intersections and on their genericity in some relevant cases is provided in Appendix~\ref{gen}.

To conclude this section, let us make some remarks on the regularity properties and the asymptotic behavior of the eigenfamilies of $H(\uu)$ in our setting. Notice that in general the regularity properties of the Hamiltonian induce similar regularity  properties of the eigenfamilies, see Proposition~\ref{proj smooth}. 
In particular, thanks to 
the Lipschitz continuity of the eigenvalues, \eqref{formcono} holds true in a neighborhood of a conical intersection, that is 
there exists a suitably small neighborhood $U$ of $\bar\uu$ and $C>0$ such that  
\bqn
\lam_{j+1}(\uu)-\lam_j(\uu)\geq C| \uu - \bar \uu|,\ \forall \uu\in U.
\label{tesi}
\eqn
Moreover, it is well known that the eigenvectors can be chosen analytic along straight lines $\uu(\cdot)$ possibly passing through eigenvalues intersections 
(see~\cite{katino},\cite[Theorem XII.13]{reed_simon}).

Consider 
a $\con^1$ curve $\uu: I \to \mathbb{R}^3$ and assume that the eigenvalues
$\Lambda_l:I \to \mathbb{R}$
and the eigenstates $\Phi_l:I \to \mathcal{H} , \ l=0,\ldots,k$  are  $\con^1(I)$. 
By direct computations we obtain that for all $t\in I$ the following equations hold: 
\begin{gather}
\dot{\Lambda}_l(t) = \langle \Phi_l(t), \big(\dot{u}_1(t) H_1+\dot{u}_2(t) H_2+\dot{u}_3(t) H_3 \big) \Phi_l(t) \rangle \label{obs1}
\\
(\Lambda_m(t)-\Lambda_l(t)) \:\langle \Phi_l(t),\dot{\Phi}_m(t)\rangle
=\langle \Phi_l(t), \big(\dot{u}_1(t) H_1+\dot{u}_2(t) H_2+\dot{u}_3(t) H_3 \Big) \Phi_m(t) \rangle.\label{obs2}
\end{gather}
If $\Lambda_j(\bar \uu)=\Lambda_{j+1}(\bar \uu)$, then, thanks to~\eqref{obs2}, for every half-line $r_{\vv}(t)=\bar \uu +t \vv$ 
 with $\vv=(v_1,v_2,v_3)$  unit vector and $t\geq 0$, we have
\beq \label{limitzero}
\lim_{t\rightarrow 0^+} \langle \Phi_j(r_{\vv}(t)), \big( v_1 H_1+v_2 H_2+v_3 H_3 \big)  \Phi_{j+1}(r_{\vv}(t)) \rangle=0.
\eeq

\subsection{The adiabatic theorem} \label{subs ad}

In this section we recall a classical formulation of the time-adiabatic theorem (\cite{bofo,kato,nenciu,panati}) adapted  to our framework. 
For a general overview see the monograph~\cite{teufel}. 

Let $H(\uu)= H_0 + \sum_{i=1}^3 u_i H_i$
satisfy {\bf (H0)}-{\bf (H1)}.
Assume that the map $I=[\tau_0,\tau_f] \ni \tau \mapsto \uu(\tau)=(u_1(\tau),u_2(\tau),u_3(\tau))$ belongs to $\con^2(I,\R^3)$.
We introduce a small parameter $\eps>0$ that controls the time scale, and the slow Hamiltonian $H(\uu(\eps t)),$  $t \in [\tau_0/\eps,\tau_f/\eps]$.
In this notations, $\tau$ is a geometric parameter used to describe the curve in the space of controls, while $t$ is the ``chronological'' time of the 
evolution along the control path.

We denote by $U^{\eps}(t,t_0)$ the time evolution (from $t_0=\tau_0/\eps$ to $t=\tau/\eps$) generated by $H(\uu(\eps t))$, 
and with $U^{\eps}_a(t,t_0)$ the time evolution generated by
the Hamiltonian $H_a(\eps t)$, where $H_a (\tau)= H(\uu(\tau)) -i\eps P_*(\uu(\tau)) \dot{P}_*(\uu(\tau)) -i \eps P^{\bot}_*(\uu(\tau)) \dot{P}^{\bot}_*(\uu(\tau))$ is the \emph{adiabatic Hamiltonian},  
$P_*(\uu)$ denotes the spectral projection of $H(\uu)$ on $\Sigma(\uu)$,  and  $P^{\bot}_*(\uu)=\mathrm{id}-P_*(\uu)$.

\bt \label{th: adiabatic}
Assume that
$%\widehat
{H}(\uu)= H_0 + \sum_{i=1}^3 u_i H_i$ satisfies {\bf (H0)}-{\bf (H1)}. %, and that $\Sigma$ is a \mariospectrum\ on $\omega \subset \mathbb{R}^3$. 
Let $I \subset \mathbb{R}$ and 
$\uu : I \rightarrow \omega$ be a $\con^2$ curve.
Then $P_* \in \con^2(I,\mathcal{L}(\mathcal{H}))$ and there exists a constant $C>0$ such that for all $\tau_0,\tau\in I$, and setting $t_0=\tau_0/\eps,t=\tau/\eps,$
\begin{equation} \label{eq: adiabatic}
\|U^{\eps}(t,t_0)-U_a^{\eps}(t,t_0)\| \leq C\eps \left(1+ \eps|t-t_0|\right).
\end{equation}
\et

\brem \label{nenciu}
If there are more than two parts of the spectrum which are separated by a gap, then it is possible to generalize the adiabatic Hamiltonian as % in the following way 
(\cite{nenciu})
$ H_a(\tau) = H(\uu(\tau)) -i \eps \sum_{\alpha} P_{\alpha}(\uu(\tau)) \dot{P}_{\alpha}(\uu(\tau))$, 
where each $P_{\alpha}(\uu(\tau))$ is the spectral projection associated with the separated portion of the spectrum labeled by $\alpha$. 
\erem

Let us now assume that $\Sigma=\{\lambda_j,\lambda_{j+1}\}$; we can take advantage of the adiabatic theorem to decouple the dynamics associated with the band  $\Sigma$
from those associated with the rest of the spectrum, in order to focus on the former. 
 
Let $\Hdue(\tau)$ denote the subspace spanned by the eigenstates associated with  $\lambda_j(\uu(\tau))$ and $\lambda_{j+1}(\uu(\tau))$. 
Since  $\Hdue(\tau)$ is two-dimensional for any $\tau$, it is possible to map it isomorphically on $\mathbb{C}^2$ and identify an effective Hamiltonian whose 
evolution is a representation of $U_a^{\eps}(t,t_0)|_{\Hdue(\eps t_0)}$ on $\mathbb{C}^2$.
In particular, if we can find a $\con^1$ eigenstate basis $\{\Phi_{1}(\uu(\tau)),\Phi_{2}(\uu(\tau))\}$ of $\Hdue(\tau)$ 
(associated with a reordering $\{\Lambda_{1}(\uu(\tau)),\Lambda_{2}(\uu(\tau))\}$ of $\{\lambda_{j}(\uu(\tau)),\lambda_{j+1}(\uu(\tau))\}$), 
then the isomorphism $\mathcal{U}(\tau) : \Hdue(\tau) \rightarrow \mathbb{C}^2$ is continuous.
Represented in $\mathbb{C}^2$,
the evolution  $U_a^{\eps}(\cdot,t_0)|_{\Hdue(\eps t_0)}$ %inside $ \Hdue(\eps t)$   
is governed by the Hamiltonian $\Heff (\eps t)$, where $\Heff (\cdot)$ is the \emph{effective Hamiltonian}, whose form is
\begin{equation}
\Heff (\tau)=	
\begin{pmatrix} \label{Heff}
\Lambda_{1}(\uu(\tau)) & 0\\
 0 & \Lambda_{2}(\uu(\tau))
\end{pmatrix}
- i \eps 
\begin{pmatrix}
\langle\Phi_{1}(\uu(\tau)), \dot{\Phi}_{1}(\uu(\tau))\rangle   & \!\! \langle\Phi_{2}(\uu(\tau)),\dot{\Phi}_{1}(\uu(\tau))\rangle\\
\langle\Phi_{1}(\uu(\tau)), \dot{\Phi}_{2}(\uu(\tau))\rangle   & \!\! \langle\Phi_{2}(\uu(\tau)), \dot{\Phi}_{2}(\uu(\tau))\rangle
\end{pmatrix},
\end{equation}
with associated propagator
$\Ueff(t,t_0) = \mathcal{U}(\eps t)U_a^{\eps}(t,t_0)\mathcal{U}^*(\eps t_0).$ 

 Theorem~\ref{th: adiabatic} implies the following.
\bt \label{effective}
Assume that $\{\lambda_j,\lambda_{j+1}\}$ is a \mariospectrum\ on some $\omega\in 
\mathbb{R}^3$ and let $\uu : I \rightarrow \omega$ be a $\con^2$ curve 
such that there exists a $\con^1$-varying  basis of $\Hdue(\cdot)$ made of eigenstates of $H(\uu(\cdot))$.
Then there exists a constant $C$ such that for all  $\tau_0,\tau\in I$, and setting $t_0=\tau_0/\eps,t=\tau/\eps,$
\[
\|\left( U^{\eps}(t,t_0)- \mathcal{U}^*(\eps t) \Ueff (t,t_0) \mathcal{U}(\eps t_0)\right)|_{\Hdue(\eps t_0)} \|
 \leq C \eps (1+\eps|t-t_0|). \]
\et

\section{Conical Intersections and general control strategy}
\label{s-conic}

\subsection{%Geometric 
Properties of conical intersections}

\medskip
Conical intersections have a characterization in terms of the non-degeneracy of a particular matrix, which contains some geometric properties of the eigenspaces relative to 
the intersecting eigenvalues, as shown below.%the following results show. 

\bdeff
We define the \emph{\mariomatrix} associated with two orthonormal elements $\psi_1,\psi_2\in \mathcal{D}(H_0)$ 
as
\[
\mathcal{M}(\psi_1,\psi_2)=\begin{pmatrix}
\langle \psi_1, H_1 \psi_2 \rangle &  \langle \psi_1, H_1 \psi_2 \rangle^* &  \langle \psi_2, H_1 \psi_2 \rangle-\langle \psi_1, H_1 \psi_1 \rangle \\ 
\langle \psi_1, H_2 \psi_2 \rangle & \langle \psi_1, H_2 \psi_2 \rangle^* &  \langle \psi_2, H_2 \psi_2 \rangle-\langle \psi_1, H_2 \psi_1 \rangle \\
\langle \psi_1, H_3 \psi_2 \rangle & \langle \psi_1, H_3 \psi_2 \rangle^* &  \langle \psi_2, H_3 \psi_2 \rangle-\langle \psi_1, H_3 \psi_1 \rangle
\end{pmatrix} .
\]

\edeff

\begin{lemma}
The quantity $\det \mathcal{M}(\psi_1,\psi_2)$ is purely imaginary and the function $(\psi_1,\psi_2) \mapsto \det{\cal M}(\psi_1,\psi_2)$ is invariant under unitary  
transformation of the argument, that is if $( \widehat{\psi}_1, \widehat{\psi}_2)^T= \mathbb{U}(\psi_1, \psi_2)^T$ for a pair $\psi_1,\psi_2$ of orthonormal elements of 
$\mathcal{D}(H_0)$ %$\mathcal{H}$ 
and $\mathbb{U}\in \textsf{U(2)} $, then 
one has $\det \mathcal{M}(\widehat{\psi}_1,\widehat{\psi}_2)=\det \mathcal{M}(\psi_1,\psi_2)$.

\label{tuttouguale}
\end{lemma}

\proof
The fact that the determinant is purely imaginary comes from direct computations. To prove its invariance, we set 
\begin{equation} \label{matrice su2}
\mathbb{U}=\begin{pmatrix} e^{i(\beta_1+\beta_3)}\cos\al &  e^{i(\beta_2+\beta_3)}\sin\al \\ -e^{-i\beta_2}\sin\al & e^{-i\beta_1} \cos\al\end{pmatrix}
\end{equation}
for some real scalars $\beta_1,\beta_2,\beta_3$, so that $\widehat{\psi}_1=e^{i\beta_3}\left(e^{i\beta_1}\cos\al\: \psi_1 + e^{i\beta_2}\sin\al \:\psi_2\right)$ and
$\widehat{\psi}_2=-e^{-i\beta_2}\sin\al \:\psi_1 + e^{-i\beta_1}\cos\al\: \psi_2$.

By direct computations it follows that
\[
\mathcal{M}(\widehat{\psi}_1,\widehat{\psi}_2)= \mathcal{M}(\psi_1,\psi_2) 
\begin{pmatrix} 
e^{-2i\beta_1-i\beta_3} \cos^2 \alpha & -e^{2i\beta_2+i\beta_3} \sin^2 \alpha & -e^{i(\beta_2-\beta_1)} \sin 2\alpha\\ 
-e^{-2i\beta_2-i\beta_3} \sin^2 \alpha & e^{2i\beta_1+i\beta_3} \cos^2 \alpha & -e^{i(\beta_1-\beta_2)} \sin 2\alpha\\ 
e^{-i(\beta_1+\beta_2+\beta_3)} \cos \alpha \sin \alpha & e^{i(\beta_1+\beta_2+\beta_3)} \cos \alpha \sin \alpha &  \cos 2\alpha
\end{pmatrix},
\]
where the second matrix on the right-hand side of the equation above has determinant equal to one.
\eproof 

\medskip

As a consequence of the result here above, the determinant of $\mathcal{M}(\psi_1,\psi_2)$ depends only on the complex space spanned by $\psi_1$ and 
$\psi_2$. 
Therefore, in a neighborhood of a conical intersection between the levels $\lambda_j,\lambda_{j+1}$  we can define the following function:
\begin{equation}
F(\uu)=\det \mathcal{M}(\psi_1,\psi_2)
\end{equation}
where $\{\psi_1,\psi_2\}$ is an orthonormal basis for the sum of eigenspaces relative to the two crossing levels. In particular, outside the intersection we can take, for instance, $\psi_1=\phi_j$ and $\psi_2=\phi_{j+1}$. 

If the levels $\lambda_j,\lambda_{j+1}$ are (locally) separated from the rest of the spectrum, 
%By the continuity of 
the projection associated with the sum of the eigenspaces of the intersecting levels is continuous with respect to $\uu$ (see Proposition~\ref{pdiff}), 
which implies that $F$ is continuous (see~\cite{adia-nostro}). 

The following result characterizes conical intersections in terms of the \mariomatrix.

\bp 
\label{iff conical}
Assume that $\{\lam_j, \lam_{j+1}\}$ is a \mariospectrum\ with $\lam_j(\bar \uu) = \lam_{j+1}(\bar \uu)$.
Let $\{\psi_1,\psi_2\}$ be an orthonormal basis of the eigenspace associated with the double eigenvalue. Then $\bar\uu$ is a conical intersection if and only if ${\cal M}(\psi_1,\psi_2)$ is nonsingular.
\ep

\proof
Let $r_{\vv}(t)=\bar\uu+t \vv$, where $\vv$ is a unit vector in $\mathbb{R}^3$, 
and let $\phi_j^{\vv},\phi_{j+1}^{\vv}$ be the limits of $\phi_j(r_{\vv}(t)),\phi_{j+1}(r_{\vv}(t))$ as $t\rightarrow 0^+$ (recall that the eigenfunctions $\phi_j,\phi_{j+1}$ can be chosen analytic along $r_{\vv}$ for $t\geq 0$). 
Assume that the intersection is not conical.
%, that is, thanks to equation~\eqref{conical equiv},
Then for every $\eps>0$ there is a unit vector $\vv^{\eps}=(v^{\eps}_1,v^{\eps}_2,v^{\eps}_3)$
such that
\[
\frac{d}{dt} \Big|_{t=0^+}\Big[\lam_{j+1}(r_{\vv_{\eps}}(t))-\lam_{j}(r_{\vv_{\eps}}(t)) \Big]  \leq \eps,
\]
that is
\[a^{\vv_\eps}:=\sum_{i=1}^3
v_i^{\eps} \left( \langle \phi^{\vv_\eps}_{j+1}, H_i \phi^{\vv_\eps}_{j+1}\rangle -\langle \phi^{\vv_\eps}_{j}, H_i \phi^{\vv_\eps}_{j}\rangle \right) 
\leq \eps,\]
while~\eqref{limitzero} implies that 
$b^{\vv_\eps}:=\sum_{i=1}^3 v_i^{\eps} \langle \phi^{\vv_\eps}_{j}, H_i \phi^{\vv_\eps}_{j+1}\rangle
=0$.
Consider an orthogonal matrix  $\mathbb{A}_{\eps}$ having  $\vv_\eps$ as first row.
Since $|\det \mathbb{A}_{\eps}|=1$, we have that
\begin{align*}
\left| \det \mathcal{M}(\phi^{\vv_\eps}_{j},\phi^{\vv_\eps}_{j+1}) \right| &=  
\left| \det \left( \mathbb{A}_{\eps}\mathcal{M}(\phi^{\vv_\eps}_{j},\phi^{\vv_\eps}_{j+1})\right) \right|= \left| \det 
\begin{pmatrix}
0 & 0 & a^{\vv_\eps}\\
c_1 & c_1^* & d_1 \\ 
c_2 & c_2^* & d_2
\end{pmatrix}\right|
=\left| 2i a^{\vv_\eps} \mathfrak{Im}(c_1c_2^*)\right|\\
&\leq C\eps (\alpha |\lambda_j(\uu)|+\beta)^2, 
\end{align*}
with $C$, $\alpha$ and $\beta$ suitable positive constants, where we have used the fact that 
\[|c_i|\leq |( \langle \phi^{\vv_\eps}_{j}, H_1\phi^{\vv_\eps}_{j+1}\rangle, \langle \phi^{\vv_\eps}_{j}, H_2 \phi^{\vv_\eps}_{j+1}\rangle, \langle \phi^{\vv_\eps}_{j}, H_3 \phi^{\vv_\eps}_{j+1}\rangle)|\] 
for $i=1,2$ and that $H_i$ is $H(\uu)$-bounded. Thus arbitrariness of $\eps$ and Lemma~\ref{tuttouguale} imply that ${\cal M}(\psi_1,\psi_{2}) $ is singular for any orthonormal basis $\{\psi_1,\psi_{2}\}$ of the double eigenspace.

Let us now prove the converse statement: assume that $\bar \uu$ is a conical intersection and, by contradiction, that ${\cal M}(\psi_1,\psi_2)$ is singular for every orthonormal basis of the eigenspace associated with the double eigenvalue.
We introduce the matrix
\[\widetilde{\mathcal{M}}(\psi_1,\psi_2)=
\begin{pmatrix}
\mathfrak{Re}\left(\langle \psi_1, H_1 \psi_2 \rangle \right) & \mathfrak{Im}\left( \langle \psi_1, H_1 \psi_2 \rangle  \right)& \left( \langle \psi_2, H_1 \psi_2 \rangle-\langle \psi_1, H_1 \psi_1 \rangle \right)\\ 
\mathfrak{Re}\left(\langle \psi_1, H_2 \psi_2 \rangle \right) &\mathfrak{Im}\left( \langle \psi_1, H_2 \psi_2 \rangle \right) &\left(  \langle \psi_2, H_2 \psi_2 \rangle-\langle \psi_1, H_2 \psi_1 \rangle \right)\\
\mathfrak{Re}\left(\langle \psi_1, H_3 \psi_2 \rangle \right) &\mathfrak{Im}\left( \langle \psi_1, H_3 \psi_2 \rangle \right)& \left( \langle \psi_2, H_3 \psi_2 \rangle-\langle \psi_1, H_3 \psi_1 \rangle \right)
\end{pmatrix},
\]
and we notice that $\det \mathcal{M}(\psi_1,\psi_2)=-2 i\det\widetilde{\mathcal{M}}(\psi_1,\psi_2)$ so that $\mathcal{M}(\psi_1,\psi_2)$ is singular if and only if $\widetilde{\mathcal{M}}(\psi_1,\psi_2)$ is. 	
The condition of conical intersection implies that %the third column of ${\cal M}(\phi^{\vv}_{j},\phi^{\vv}_{j+1})$ (and of $\widetilde{\mathcal{M}}(\phi^{\vv}_{j},\phi^{\vv}_{j+1})$) is non-null, and that 
%\begin{equation} \label{conicalbis}
\[
a^{\vv}=\sum_{i=1}^3
v_i \left( \langle \phi^{\vv}_{j+1}, H_i \phi^{\vv}_{j+1}\rangle -\langle \phi^{\vv}_{j}, H_i \phi^{\vv}_{j}\rangle \right)  \neq 0,\quad\forall \vv\in \mathbb{R}^3.
\]
%\end{equation}
Moreover, equation~\eqref{limitzero} implies that 
$
b^{\vv}=\sum_{i=1}^3 v_i \langle \phi^{\vv}_{j}, H_i \phi^{\vv}_{j+1}\rangle
=0$,
so that the third column of the matrix $\widetilde{\mathcal{M}}(\phi^{\vv}_{j},\phi^{\vv}_{j+1})$ is never linearly dependent from the first two.
In particular the matrix $\widetilde{\mathcal{M}}(\phi^{\vv}_{j},\phi^{\vv}_{j+1})$ is singular only if the first two columns of the matrix are linearly dependent. 
%which is equivalent to say that $\langle \phi^{\vv}_{j}, H_i \phi^{\vv}_{j+1}\rangle$ is either a real for every $i$ or purely imaginary for every $i$. 
Thus, up to multiplying $\phi^{\vv}_{j}, \phi^{\vv}_{j+1}$ by a phase factor, we can always assume that $\langle \phi^{\vv}_{j}, H_i \phi^{\vv}_{j+1}\rangle\in \mathbb{R}, \ i=1,2,3$. %, and that the second column of $\widetilde{\mathcal{M}}(\phi^{\vv}_{j},\phi^{\vv}_{j+1})$ is null.

Let us now fix $\vv \in \R^3$ and let us call $W$ the orthogonal complement in $\mathbb{R}^3$ of the vector $\left(\langle \phi^{\vv}_{j}, H_1 \phi^{\vv}_{j+1}\rangle,\langle \phi^{\vv}_{j}, H_2 \phi^{\vv}_{j+1}\rangle,\langle \phi^{\vv}_{j}, H_3 \phi^{\vv}_{j+1}\rangle \right)$. We have that $\vv \in W$ and $\dim W\geq2$.
It is easy to prove that for every $\ww \in W$ the limit basis 
$\{\phi^{\ww}_{j}, \phi^{\ww}_{j+1}\}$ is equal to $\{\phi^{\vv}_{j}, \phi^{\vv}_{j+1}\}$, up to exchanges between the two elements and up to phases. Indeed, by definition of $W$ we have that 
$\langle \phi^{\vv}_{j},(w_1H_1+w_2 H_2 + w_3 H_3) \phi^{\vv}_{j+1}\rangle=0$, and  writing $(\phi^{\vv}_{j}, \phi^{\vv}_{j+1})^T=\mathbb{U}(\phi^{\ww}_{j}, \phi^{\ww}_{j+1})^T$ for some $\mathbb{U}$ of the form~\eqref{matrice su2}, we obtain that 
\begin{align*}
0&=\langle \phi^{\vv}_{j},(w_1H_1+w_2 H_2 + w_3 H_3) \phi^{\vv}_{j+1}\rangle \\
&=e^{-i\beta_3}\left( (e^{-2i\beta_1}\cos^2 \alpha-e^{-2i\beta_2} \sin^2 \alpha)\ b^{\ww} +  e^{-i(\beta_1+\beta_2)} \sin \alpha \cos \alpha \ a^{\ww} \right), 
\end{align*}
therefore it must be $\alpha=k\frac{\pi}{2},\ k\in \mathbb{Z}$, since $b^{\ww}=0$ and $a^{\ww}\neq 0$.

Let us now consider the vector 
\[
\Upsilon= \left(\langle \phi^{\ww}_{j+1},H_1 \phi^{\ww}_{j+1}\rangle-\langle \phi^{\ww}_{j},H_1 \phi^{\ww}_{j}\rangle,\langle \phi^{\ww}_{j+1},H_2  \phi^{\ww}_{j+1}\rangle -  \langle \phi^{\ww}_{j},H_2 \phi^{\ww}_{j}\rangle,\langle \phi^{\ww}_{j+1},H_3 \phi^{\ww}_{j+1}\rangle -\langle \phi^{\ww}_{j},H_3 \phi^{\ww}_{j}\rangle \right),
\]
which corresponds to the third column of $\widetilde{\mathcal{M}}(\phi^{\ww}_{j},\phi^{\ww}_{j+1})$.
Notice that the definition of $\Upsilon$, up to a sign, does not depend on the choice of $\ww\in W$, by what precedes. 
Since $\Upsilon^{\bot}$, the orthogonal complement of $\Upsilon$ in $\mathbb{R}^3$, has dimension 2  
there exists a non-zero $\widetilde{\ww} \in W \cap \Upsilon^{\bot}$.
By definition of $\widetilde{\ww}$ 
we have that $a^{\widetilde{\ww}}=\widetilde{\ww}\cdot \Upsilon=0$.

We get a contradiction, thus the matrix $\widetilde{\mathcal{M}}(\phi^{\vv}_{j},\phi^{\vv}_{j+1})$ must be nonsingular, and therefore
also $\mathcal{M}(\phi^{\vv}_{j},\phi^{\vv}_{j+1})$ has to.
\eproof

\medskip

As an example, consider the Hamiltonian~$H(\uu)\in i\mathfrak{u}(3)$, with
\begin{equation}
H_0=
\left(
\begin{array}{ccc}
0 &0 &0 \\
0 & 0 & 0 \\
0& 0 & 1
\end{array}
\right)
,\ 
H_1=
\left(
\begin{array}{ccc}
1 & i  & 0 \\
-i  & 0 & 1 \\
0 & 1 & -1
\end{array}
\right)
,\ 
H_2=
\left(
\begin{array}{ccc}
0 & 0  & i \\
0 & 1 & 0 \\
-i & 0 & 0
\end{array}
\right)
,\ 
H_3=
\left(
\begin{array}{ccc}
-1 & 1 & -1 \\
1  & 1 & 0 \\
-1 & 0 & 0
\end{array}
\right).\label{H-ricci}
\end{equation}
The Hamiltonian $H(\cdot)$ admits a double eigenvalue at $\uu=0$ corresponding to the two lowest levels. 
A simple computation leads to $\det \mathcal{M}(e_1,e_2)=-2i$ where $e_1=(1,0,0)^T,e_2=(0,1,0)^T$ form a basis of the double eigenspace at $0$. 
Thus the eigenvalue intersection is conical.
An example of conical intersection in an infinite dimensional setting is provided in Appendix~\ref{s-infinity}.

\medskip

%As noticed above,

A peculiarity  of conical intersections is that, when approaching the singularity from different directions, the eigenstates corresponding to the intersecting eigenvalues have different limits. The following proposition provides the relation between these limits.

\bp \label{change of basis}
Let $\bar \uu$ be a conical intersection between $\lambda_j$ and $\lambda_{j+1}$.
Let ${\bf v}_0,\vv\in \mathbb{R}^3$ be two unit vectors, and call
$\phi_j^0,\phi_{j+1}^0$ the limits as $t \to 0^+$ of the eigenstates $\phi_j(r_0(t) ),\phi_{j+1}(r_0(t) )$ along a straight line $r_0(t)=\bar \uu+t{\bf v}_0$,
and $\phi_j^{\bf v},\phi_{j+1}^{\bf v}$ the limit basis  along the straight line $r_{\bf v}(t)=\bar \uu+t{\bf v}$.
Then, up to phases, the following relation holds:
\begin{equation} \label{tranx}
\begin{pmatrix}
\phi_j^{\bf v} \\ \phi_{j+1}^{\bf v}
\end{pmatrix}=
\begin{pmatrix}
\cos \bth &  e^{-i \beta}\sin \bth \\
-e^{i \beta}\sin \bth  & \cos \bth  
\end{pmatrix}
\begin{pmatrix}
\phi_j^0 \\ \phi_{j+1}^0
\end{pmatrix},
\end{equation}
where the parameters $\bth=\bth({\bf v} )$ and $\beta=\beta({\bf v} )$ satisfy the following equations:
\begin{gather}
\tan 2 \bth({\bf v})  = (-1)^\xi \frac{ 2 |\langle \phi_j^0, H_{\bf v} \phi_{j+1}^0 \rangle | }{\langle \phi_j^0, H_{\bf v} \phi_{j}^0\rangle-\langle \phi_{j+1}^0, H_{\bf v} \phi_{j+1}^0\rangle}  \label{alpha} \\
\beta({\bf v} ) \stackrel{(\!\!\!\!\!\!\mod\! 2\pi\!)}{=} \arg \langle \phi_j^0, H_{\bf v} \phi_{j+1}^0\rangle+\xi\pi,\label{beta}
\end{gather}
where $H_{\bf v}= \sum_{i=1}^m H_i  v_i$ and $\xi=0,1$.
\ep

\proof
First of all, we notice that all pairs of orthonormal 
eigenstates of $H(\bar \uu)$ relative to the degenerate eigenvalue can be obtained by the action of the group \textsf{U(2)}  through the pair $(\phi_j^0,\phi_{j+1}^0)$. 
Nevertheless, we are not interested on the global phases of the states, that is we consider the equivalence relation 
$(e^{i\beta_1}\psi_1,e^{i\beta_2}\psi_2) \sim (\psi_1,\psi_2)\ \forall \; \beta_1,\beta_2\in \mathbb{R}$. 
Therefore for any unit vector $\vv\in \mathbb{R}^3$ we can obtain a representative of the pair $(\phi_j^{\bf v},\phi_{j+1}^{\bf v})$ through the transformation~\eqref{tranx}, for some $\bth(\vv)$ and some $\beta(\vv)$. This gives
\begin{align}
0 &= \langle \phi_j^{\bf v}, H_{\bf v} \phi_{j+1}^{\bf v} \rangle  \nonumber \\
&= \left(\langle \phi_j^0, H_{\bf v} \phi_{j+1}^0 \rangle,\langle \phi_j^0, H_{\bf v} \phi_{j+1}^0 \rangle^*,
\langle \phi_{j+1}^0, H_{\bf v} \phi_{j+1}^0 \rangle-\langle \phi_j^0, H_{\bf v} \phi_{j}^0 \rangle \right) 
\begin{pmatrix}
\cos^2\bth({\bf v})\\
-e^{2i\beta({\bf v})}\sin^2\bth({\bf v})\\
e^{i\beta({\bf v})}\cos\bth({\bf v})
\sin\bth({\bf v})  
\end{pmatrix} \nonumber \\
&= e^{i\beta({\bf v})} \Big(\cos^2\bth({\bf v}) |\langle\phi_j^0, H_{\bf v} \phi_{j+1}^0 \rangle| e^{i(\gamma -  \beta({\bf v}))} - \sin^2\bth({\bf v}) 
|\langle\phi_j^0, H_{\bf v} \phi_{j+1}^0 \rangle|e^{-i(\gamma -  \beta({\bf v}))} + \nonumber  \\
&+ \cos\bth({\bf v})\sin\bth({\bf v})(\langle \phi_{j+1}^0, H_{\bf v} \phi_{j+1}^0 \rangle-\langle \phi_j^0, H_{\bf v} \phi_{j}^0 \rangle) \Big),\label{ottobre rosso}
\end{align}
with $\gamma=\arg\left(\langle \phi_j^0, H_{\bf v} \phi_{j+1}^0 \rangle \right)$, whenever  $ \langle\phi_j^0, H_{\bf v} \phi_{j+1}^0 \rangle \neq 0$.

If $\vv= \pm \vv_0$, then $\langle\phi_j^0, H_{\bf v} \phi_{j+1}^0 \rangle=0$, and the conicity condition implies that $\bth(\vv)=k\pi/2$. In particular,
$\beta(\vv)$ can be any real number.

If $\vv$ is not parallel to $\vv_0$, then $\langle\phi_j^0, H_{\bf v} \phi_{j+1}^0 \rangle\neq 0$ and $\beta(\bf v)$ can only take the values $\arg \langle \phi_j^0, H_{\bf v} \phi_{j+1}^0\rangle+\xi\pi,\ \xi=0,1$. Indeed, since the term $\langle \phi_{j+1}^0, H_{\bf v} \phi_{j+1}^0 \rangle-\langle \phi_j^0, H_{\bf v} \phi_{j}^0 \rangle$ is real,
the imaginary part of the term into parenthesis of~\eqref{ottobre rosso} %shall be
is zero:
\[
|\langle\phi_j^0, H_{\bf v} \phi_{j+1}^0 \rangle| \sin(\gamma -  \beta({\bf v}))=0.
\]
If we choose $\beta({\bf v})=\gamma$, then we can prove by computation that $\bth({\bf v})$ %shall 
must
satisfy 
\begin{equation} \label{theta 1}
\tan 2 \bth({\bf v})  = \frac{ 2 |\langle \phi_j^0, H_{\bf v} \phi_{j+1}^0 \rangle | }{\langle \phi_j^0, H_{\bf v} \phi_{j}^0\rangle-\langle \phi_{j+1}^0, H_{\bf v} \phi_{j+1}^0\rangle}
\end{equation}
while if $\beta({\bf v})=\gamma+\pi$,
$\bth({\bf v})$ %shall 
must satisfy 
\begin{equation} \label{theta 2}
\tan 2 \bth({\bf v})  = -\frac{ 2 |\langle \phi_j^0, H_{\bf v} \phi_{j+1}^0 \rangle | }
{\langle \phi_j^0, H_{\bf v} \phi_{j}^0\rangle-\langle \phi_{j+1}^0, H_{\bf v} \phi_{j+1}^0\rangle}.
\end{equation}
Notice that the two pairs $(\bth({\bf v}),\beta({\bf v}))$ and $(-\bth({\bf v}),\beta({\bf v})+\pi)$ give the same transformation in~\eqref{tranx}.
\eproof 

\medskip

It can be seen that not all the solutions of~\eqref{alpha}-\eqref{beta} provide the correct transformation~\eqref{tranx}, which, 
nevertheless, is easy to detect. The good solutions of~\eqref{alpha}-\eqref{beta} constitute four branches which are continuous with respect to $\vv$, and they can be constructed as follows.
Let ${\bf w}(s),\ s\in [0,\bar s]$, be a curve joining $\bf{v}_0$ to $\vv$ such that 
${\bf w}(s) \notin \{ {\bf v}_0,-{\bf v}_0 \}$ for every $s\in (0,\bar s)$; 
for conical intersections, it is possible to associate with such a curve % every curve ${\bf w}(s)$ of unit vectors joining ${\bf v}_0$ to ${\bf v}$ 
a continuous solution $(\bth({\bf w} (s)),\beta({\bf w} (s)))$ of~\eqref{alpha}-\eqref{beta} with $\bth({\bf v}_0)=0$ and compatible with~\eqref{tranx}.
In particular, if we choose $\bth(\vv)$ according to~\eqref{theta 1}, it is easy to see that $\bth({\bf w} (s))\in [-\pi/2,0]$ for $s\in [0,\bar s]$,
 from which one deduces that the final value $\bth(\vv)=\bth({\bf w} (\bar s))$ is independent of the chosen path and continuously depends on 
$\vv$. Moreover, it is easy to see that $\bth(-{\bf v}_0)=-\pi/2$. 
Similarly, one can show that $\beta(\vv)=\beta({\bf w} (\bar s))$ is independent of the chosen path and continuous outside $\{ {\bf v}_0,-{\bf v}_0 \}$.  
Note that the fact that $\beta$ is discontinuous at $-{\bf v}_0$ implies that the corresponding limit basis $(\phi_j^{\bf v},\phi_{j+1}^{\bf v})$ 
has a discontinuity at $-{\bf v}_0$, that is, its limit depends on the path.

We can repeat the same argument choosing $\bth(\vv)$ according to~\eqref{theta 2}, with $\bth(\vv_0)=0$.
% in which case $\bth({\bf w} (s))\in [0,\pi/2]$ for $s\in [0,\bar s]$.
The other two continuous branches are obtained choosing the initial condition $\bth(\vv_0)=\pi$.

\subsection{The basic control algorithm}
\label{control-trivial}
Let us consider the following controllability problem.
\begin{quote}
%\begin{center}
Let $H(\cdot)$ satisfy {\bf (H0)}-{\bf (H1)} and $\Sigma(\cdot)=\{\lambda_0(\cdot),\dots,\lambda_k(\cdot)\}$. Then, given $\eps>0$, $\uu^{s},\uu^{f}\in\om$, $j\in\{0,\dots,k\}$ and $p \in [0,1]^{k+1}$ such that $\sum_{l=0}^k p_l^2=1$, find $T>0$ and a path $\uu:[0,T]\to\om$ with $\uu(0)=\uu^{s}$ and $\uu(T)=\uu^{f}$ such that  
\begin{equation} 
\|\psi(T)-\sum_{l=0}^k p_l e^{i\vartheta_l} \phi_l(\uu^f)\|\leq  \varepsilon,
\end{equation}
where $\psi(\cdot)$ is the solution of (\ref{hamiltonianabold}) with $\psi(0)=\phi_j(\uu^{s})$, and $\vartheta_0,\ldots,\vartheta_k \in \mathbb{R}$ are some possibly unknown phases.
%\end{center}
\end{quote}
If all  levels are connected by means of conical intersection occurring at different values of the control, the results obtained in the previous section provide the basic elements in order to construct a family of control paths solving the problem here above. This can be done by taking advantage of 
the following proposition, which describes the spreading of occupation probabilities induced when a path in the space of controls passes through a conical intersection.

\bp \label{controllabilita}
Let $\bar\uu$ be a conical intersection between the eigenvalues $\lam_j,\lam_{j+1}$. Consider the curve  $\gamma:[0,1]\rightarrow \omega$ 
defined by 
\[
\gamma(t)=
\begin{cases}
\bar\uu +(t-\tau_0) \ww_0 & t\in [0,\tau_0]\\
\bar\uu +(t-\tau_0) \vv & t\in [\tau_0,1]
\end{cases},
\]
for some 
$\tau_0\in(0,1)$ and
 some unit vectors $\ww_0,\vv$.
Then there exists $C>0$ such that, for any  $\varepsilon >0$,
\bqn 
\|\psi(1/\eps)- \pi_1e^{i\vartheta_1}\phi_j(\gamma(1))-\pi_2e^{i\vartheta_2}\phi_{j+1}(\gamma(1))\|\leq C \sqrt{\eps}
\eqn 
where $\vartheta_1,\ \vartheta_2\in \mathbb{R}$, $\psi(\cdot)$ is the solution of equation~\eqref{hamiltonianabold} with $\psi(0)=\phi_j(\gamma(0))$ corresponding to the control $\uu:[0,1/\eps]\rightarrow \omega$ defined by $\uu(t)=\gamma(\eps t)$, 
\[
\pi_1=|\cos \left( \bth({\bf v})  \right) |, \ 
\pi_2=|\sin \left( \bth({\bf v})  \right) |,
\]
and $\bth(\cdot)$ is the only solution of equation~\eqref{theta 1} such that $\bth(\vv) \in (-\pi/2,0]$ for $\vv \neq -\ww_0$, and $\bth(-\ww_0)=-\pi/2$, where the limit basis in~\eqref{theta 1}
is given by the limits 
$\phi_j(\gamma(\tau_0^-) ),\phi_{j+1}(\gamma(\tau_0^-) )$, respectively. 
\ep

\proof
We consider the Hamiltonian $H(\uu(t)), \ t \in [0,1/\eps]$.
Since the control function $\uu(\cdot)$ is not $\con^2$ at the singularity, we cannot directly apply the adiabatic theorem. Instead, we consider separately the evolution on the two subintervals (in time $t$) $[0,\tau_0 /\eps]$ and $[ \tau_0 /\eps,1/\eps]$. 

Since the eigenstates $\phi_j(\uu(t)),\phi_{j+1}(\uu(t))$ are piecewise $\con^1$, we can apply~\cite[Corollary 2.5]{teufel} 
and obtain that there exists a phase $\theta_1$ 
(depending on $\eps$) such that
\[
\|\psi(\tau_0/\eps)-e^{i\theta_1}\phi_j(\gamma(\tau_0^-))\| \leq C \sqrt{\eps},
\]
for some constant $C>0$.
By Proposition~\ref{change of basis}, this implies that 
\[
\|\psi(\tau_0/\eps)-e^{i\theta_1}\left(\cos \bth(\vv)\ \phi_j(\gamma(\tau_0^+)) -e^{-i\beta(\vv)} \sin \bth(\vv)\ \phi_{j+1}(\gamma(\tau_0^+))\right)\| \leq C \sqrt{\eps},
\]
with $\bth(\vv)$ as in the statement of the proposition and $\beta(\vv)= \arg \langle \phi_j^0, H_{\bf v} \phi_{j+1}^0\rangle$.

Applying~\cite[Corollary 2.5]{teufel} also in the time interval $(\tau_0/\eps,1/\eps]$, 
 we conclude that there exists two phases 
$\vartheta_1$ and $\vartheta_2$
(depending on $\eps$) such that
\[
\|\psi(1/\eps)-\left(e^{i\vartheta_1}\cos \bth(\vv) \ \phi_j(\gamma(1)) -e^{i\vartheta_2} \sin \bth(\vv) \ \phi_{j+1}(\gamma(1))\right)\| \leq C' \sqrt{\eps},
\]
for some constant $C'>0$.
\eproof 

\medskip 

For control purposes, it is interesting to consider the case in which the initial probability is concentrated in the first level, the final occupation probabilities $p_1^2$ and $p_2^2$ are prescribed, and we want to determine a path that induces the desired transition. For a given line reaching 
the conical intersection, the outward directions that provide the required spreading of probability are given in the following proposition. The proof follows from simple computations and is thus omitted.

\bp \label{intercono}
Let $\bar \uu$ be a conical intersection between the eigenvalues $\lambda_{j},\lambda_{j+1}$, and let $\pi_1,\pi_2$ be positive constants such that %two probability weights with 
$\pi_1^2+\pi_2^2=1$.
Consider the line $r(t)=\bar \uu + (t-\tau_0) \ww_0, \ t \in [0,\tau_0]$, for a unit vector $\ww_0 \in \R^3$  and some $\tau_0 \in (0,1)$, 
and set $\phi^0_j=\phi_j(\gamma(\tau_0^-) )$ and $\phi_{j+1}^0=\phi_{j+1}(\gamma(\tau_0^-) )$. 

Then the locus formed by the directions $\vv \in \R^3$ that give rise to transformation~\eqref{tranx} with $\pi_1=|\cos \Xi (\vv)|$ and $\pi_{2}=|\sin \Xi(\vv) |$ is
given by the following expression whenever $\pi_1\notin \{0,1\}$
\begin{equation} \label{ultracono}
    \widetilde{\mathcal{M}} (\phi_{j}^0,\phi_{j+1}^0)^{-T} (\mathcal{K}),
\end{equation}
where
$\mathcal{K}= \big\{ (x,y,z) \in \R^3 :  \sqrt{x^2+y^2}=C z \big\}$ and 
\[
C=\frac{\pi_1 \pi_{2}}{\pi_1^2-\pi_{2}^2}.
\]
Otherwise, if $\pi_1=0$ then $\vv=\ww_0$ and if $\pi_1=1$ then $\vv=-\ww_0$.
\ep

\noindent
The controllability problem presented at the beginning of this section can be solved taking advantage of the results shown above.
The strategy consists in constructing a piecewise $\con^2$ path joining $\uu^s$ with $\uu^f$ that passes through the conical intersections $\bar \uu_j$  between the $j$-th and the $(j+1)$-th levels, $j=0,\ldots,k-1$, and avoids any other degeneracy point.
 The tangent directions at the conical intersection are chosen according to the 
probability weights $p_i^2$, as explained in Proposition~\ref{intercono}.

When we are far from all the conical intersections, we 
approximate the evolution with that of the adiabatic Hamiltonian %\[
\begin{equation} \label{adiatutto}
h_a(\tau) = H(\gamma(\tau)) -i\eps \sum_{l=0}^k  P_l(\tau) \dot{P}_l(\tau) -i\eps  P^{\bot}(\tau) \dot{P}^{\bot}(\tau),
\end{equation} 
where $P_l(\tau) $ is the spectral projection onto the eigenspace relative to $\lambda_l(\gamma(\tau))$
and $P^{\bot}(\tau)=\mathrm{id} - \sum_{l=0}^k  P_l(\tau) $. %The time-derivatives above ( $\dot{ }$ ) have to be intended with respect to $\tau$.
The evolution associated with~\eqref{adiatutto} conserves the occupation probabilities relative to each energy level $\lambda_l,\ l=0,\ldots,k$, then 
so does the evolution of $H(\gamma(\tau))$, with an approximation of order $\eps$ (see Remark~\ref{nenciu}). 

In a neighborhood of the conical intersection between $\lambda_j$ and $\lambda_{j+1}$, $j=0,\ldots,k-1$, 
we decouple the evolution inside the band of the two intersecting levels from that relative to the rest of the spectrum, that is we approximate the dynamics with the 
ones associated with 
the adiabatic Hamiltonian
\begin{equation} 
h_a(\tau)  = H(\gamma(\tau)) -i\eps \boldsymbol{P}(\tau) \dot{\boldsymbol{P}}(\tau)  -i\eps \sum_{\substack{l=0 \\ l\neq j,j+1} }^k  P_l(\tau) \dot{P}_l(\tau)
-i\eps  P^{\bot}(\tau) \dot{P}^{\bot}(\tau) \label{adiaj}
\end{equation}
where $\boldsymbol{P}(\tau)$ is the spectral projection relative to $\{\lambda_j(\gamma(\tau)),\lambda_{j+1}(\gamma(\tau))\}$. 
The evolution associated with~\eqref{adiaj} conserves the occupation probability relative to the band
$\{\lambda_j,\lambda_{j+1}\}$, the ones relative to all other energy levels in $\{\lambda_0,\ldots,\lambda_k\}$ and the one associated with the rest of the spectrum. 
The evolution given by~\eqref{adiaj} inside the band
$\{\lambda_j,\lambda_{j+1}\}$ is described in~\cite[Corollary 2.5]{teufel} and in Proposition~\ref{controllabilita}.

Let us now explicitly determine the path in a neighborhood of each conical intersection.
Without loss of generality, we assume that $\psi(0)=\phi_0(\uu^s)$; in the other cases a path can be obtained similarly. 
We construct a %continuous 
(piecewise $\mathcal{C}^2$) path $\gamma:[0,1]\to \omega$ as follows. %having the following properties. 

We set $\gamma(0)=\uu^s$ and, for some $0<\tau_0<1$ we choose $\gamma|_{[0,\tau_0]}$ in such a way that $\gamma(\tau_0)=\bar \uu_0$ and all the eigenvalues $\lambda_l(\gamma(\tau))$ are simple for every $l=0,\ldots,k$ %and $\gamma(\tau)$ 
and $\tau\in [0,\tau_0)$. 
Moreover, $\gamma(\cdot)$ is chosen tangent to a segment in a neighborhood of $\bar \uu_0$, and we call $\vv_0^-$ the tangent direction of $\gamma(\cdot)$ at $t=\tau_0^-$.
In a neighborhood of $\bar \uu_0$ and for $\tau \geq \tau_0$, $\gamma(\cdot)$, is chosen to be tangent to $\vv_0^+$ where  the outward direction $\vv_0^+$ is chosen according to  Proposition~\ref{intercono} with $\ww_0=\vv_0^-$, $\pi_1=p_0$, $\pi_2=\sqrt{1-p_0^2}$.

The rest of the path is constructed recursively. Assume that the path has been defined up to time $\tau_{j-1}$, for some $j=1,\ldots,k-1$, with $\gamma(\tau_{j-1})=\bar \uu_{j-1}$, and that an outward direction $\vv_{j-1}^+$ has been selected. %inside the suitable cone. 
Then choose $\tau_j\in(\tau_{j-1},1)$  and the smooth path
$\gamma|_{(\tau_{j-1},\tau_j]}(\cdot)$ 
joining $\bar \uu_{j-1}$ with $\bar \uu_{j}$
such that all the eigenvalues $\lambda_l(\gamma(\tau))$, are simple for every $l=0,\ldots,k$ and for $\tau\in (\tau_{j-1},\tau_j)$, in a neighborhood of 
$\bar \uu_{j-1}$ the curve is tangent to $\vv_{j-1}^+$, and  in a neighborhood of $\bar \uu_{j}$ the curve is tangent to some segment that we call $\vv_j^-$.
To choose the outward direction $\vv_j^+$ of $\gamma(\cdot)$ at $\bar \uu_{j}$ we apply Proposition~\ref{intercono} with $\ww_0=\vv_j^-$, $\pi_1=\frac{p_{j-1}}{\sqrt{\sum_{l=j-1}^k p_l^2}}$ and 
$\pi_2=\sqrt{1-\pi_1^2}$.

The last arc defined on $(\tau_{k-1},1]$ is simply constructed by joining $\bar\uu_{k-1}$ with $\uu^f$, taking care of choosing $\gamma(\cdot)$  tangent to the outward direction in a neighborhood of $\bar\uu_{k-1}$.

To avoid highly non-homogeneous parameterizations, the path $\gamma(\cdot)$ can be reparameterized by arc-length.

\begin{figure}
\begin{center}
\begin{picture}(0,0)%
\includegraphics{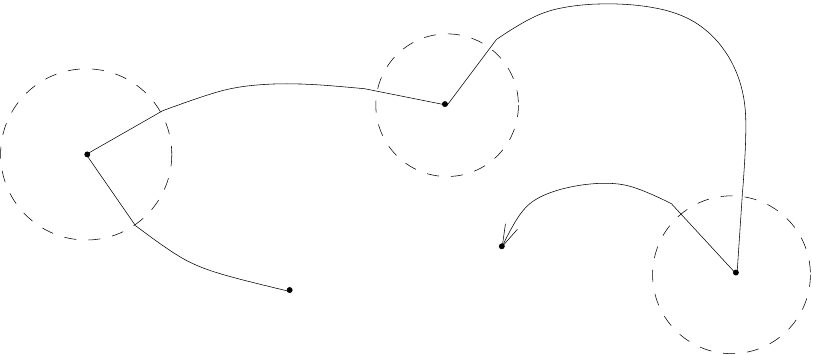}%
\end{picture}%
\setlength{\unitlength}{1381sp}%
\begingroup\makeatletter\ifx\SetFigFont\undefined%
\gdef\SetFigFont#1#2#3#4#5{%
  \reset@font\fontsize{#1}{#2pt}%
  \fontfamily{#3}\fontseries{#4}\fontshape{#5}%
  \selectfont}%
\fi\endgroup%
\begin{picture}(11122,4816)(544,-5575)
\put(1126,-2761){\makebox(0,0)[lb]{\smash{{\SetFigFont{8}{9.6}{\rmdefault}{\mddefault}{\updefault}{\color[rgb]{0,0,0}$\uu_0$}%
}}}}
\put(6391,-2506){\makebox(0,0)[lb]{\smash{{\SetFigFont{8}{9.6}{\rmdefault}{\mddefault}{\updefault}{\color[rgb]{0,0,0}$\uu_1$}%
}}}}
\put(4681,-4936){\makebox(0,0)[lb]{\smash{{\SetFigFont{8}{9.6}{\rmdefault}{\mddefault}{\updefault}{\color[rgb]{0,0,0}$\uu^s$}%
}}}}
\put(7036,-4531){\makebox(0,0)[lb]{\smash{{\SetFigFont{8}{9.6}{\rmdefault}{\mddefault}{\updefault}{\color[rgb]{0,0,0}$\uu^f$}%
}}}}
\put(10741,-4621){\makebox(0,0)[lb]{\smash{{\SetFigFont{8}{9.6}{\rmdefault}{\mddefault}{\updefault}{\color[rgb]{0,0,0}$\uu_2$}%
}}}}
\put(10876,-2236){\makebox(0,0)[lb]{\smash{{\SetFigFont{5}{6.0}{\rmdefault}{\mddefault}{\updefault}{\color[rgb]{0,0,0}$\left(\begin{array}{c} p_0\\p_1\\ \sqrt{1-p_0^2-p_1^2}\\0\end{array}\right)$}%
}}}}
\put(2551,-1261){\makebox(0,0)[lb]{\smash{{\SetFigFont{5}{6.0}{\rmdefault}{\mddefault}{\updefault}{\color[rgb]{0,0,0}$\left(\begin{array}{c} p_0\\ \sqrt{1-p_0^2}\\0\\0\end{array}\right)$}%
}}}}
\put(7951,-4111){\makebox(0,0)[lb]{\smash{{\SetFigFont{5}{6.0}{\rmdefault}{\mddefault}{\updefault}{\color[rgb]{0,0,0}$\left(\begin{array}{c} p_0\\p_1\\p_2\\p_3\end{array}\right)$}%
}}}}
\put(2401,-4936){\makebox(0,0)[lb]{\smash{{\SetFigFont{5}{6.0}{\rmdefault}{\mddefault}{\updefault}{\color[rgb]{0,0,0}$\left(\begin{array}{c} 1\\0\\0\\0\end{array}\right)$}%
}}}}
\end{picture}%
\caption{Construction of the path $\gamma(\cdot)$. The corners at the points $\uu_i$, corresponding to conical intersections, are chosen in such a way that they induce the desired spreading.}\label{path}
\end{center}
\end{figure}
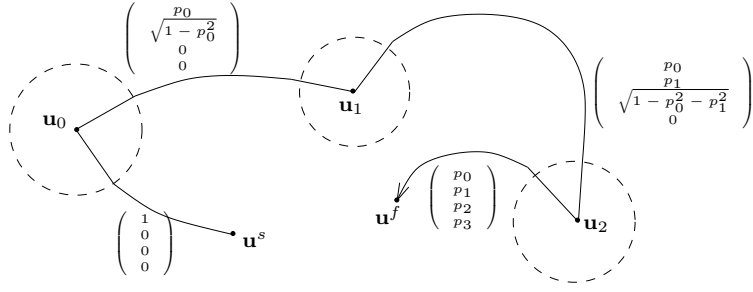

Let us now reparameterize the time setting $t=\tau/\eps$, for some small positive $\eps$. The adiabatic theorem and Proposition~\ref{controllabilita} lead to the estimate
\[\|\psi(1/\eps)-\sum_{l=0}^k p_l e^{i\vartheta_l} \phi_l(\uu^f)\|\leq C \sqrt\eps,
\]
for some $\vartheta_0,\ldots,\vartheta_k\in \mathbb{R}$ and some $C>0$ depending on the path $\gamma(\cdot)$ and on the gaps in the spectrum.
The geometric construction of the path $\gamma(\cdot)$ is represented in Figure~\ref{path}.

\section{An improvement of the efficiency of the algorithm: the non-mixing field}%\curvepao}
\label{una-a-caso}

The problem of reducing transition times is quite important in quantum control, in particular for the need of reducing decoherence. 
Therefore, in this section we will show how to do it by means of some special paths in the space of controls that eliminate, in the adiabatic approximations, the
error coming from the intersection of the eigenvalues, so that the only error induced by the application of the adiabatic theorem comes from 
the gap
among the remaining eigenvalues and those among these eigenvalues and the band of the intersecting ones.

Let us consider a pair  $\{\lam_j,\lam_{j+1}\}$ of eigenvalues, which are simple and separated from the rest of the spectrum, according to assumption \textbf{(H1)}, in a certain open set $\omega\subset \mathbb{R}^3 $, except for a point $\bar\uu \in \omega$,  where they have a conical intersection. 
We are interested in the dynamics inside the subspace $\boldsymbol{P}_{\uu}\mathcal{H}$, where $\boldsymbol{P}_{\uu}$ denotes the projection associated with the two levels $\{\lam_j(\uu),\lam_{j+1}(\uu)\}$ for $\uu\in \omega$: we know that, under adiabatic evolution, the dynamics are described by the effective Hamiltonian~\eqref{Heff}.
To improve the precision of the result, the idea is to
cancel the off-diagonal terms in the effective Hamiltonian, which are responsible of the error of order $\sqrt\eps$ the estimates given in~\cite[Corollary 2.5]{teufel} and in Proposition~\ref{intercono}. In order to do that, we choose some special trajectories  in $\omega$ along which the term $\langle \phi_{j},\dot{\phi}_{j+1}\rangle$ is null, that is, thanks to~\eqref{obs2}, we look for
curves $\gamma(\cdot)$ in the space of controls that satisfy the equation
\begin{equation} \label{paolo condition}
\langle\phi_{j}(t),(\dot{\gamma}_1(t)H_1+\dot{\gamma}_2(t) H_2+\dot{\gamma}_3(t) H_3)\phi_{j+1}(t)\rangle=0.
\end{equation}
We denote the first column of the conicity matrix $\mathcal{M}(\psi_1,\psi_2)$ by 
\[\mm(\psi_1,\psi_2)=\left(\langle \psi_1, H_1 \psi_2 \rangle,\langle \psi_1, H_2 \psi_2\rangle,\langle \psi_1, H_3 \psi_2 \rangle \right)^T,  \]
and its components $\langle \psi_1, H_i \psi_2 \rangle$ as $m_i$.
It follows by definition that the real vector 
\begin{align} \label{Xpaoprimo}
X(\psi_1,\psi_2)&=\frac{\mm(\psi_1,\psi_2) \times \mm^*(\psi_1,\psi_2)}{2i}\\
&=\left(\mathfrak{Im}(m_2 m_3^*),\mathfrak{Im}(m_3 m_1^*),\mathfrak{Im}(m_1 m_2^*)\right)^T \nonumber, 
\end{align}
where $\times$ denotes the cross product,
is orthogonal to both $\mm(\psi_1,\psi_2)$ and $\mm^*(\psi_1,\psi_2)$.

\begin{remark}
Let us remark that the vector $X(\psi_1,\psi_2)$ is invariant under phase changes in the argument, that is 
 $X(\psi_1,\psi_2)=X(e^{i\beta_1}\psi_1,e^{i\beta_2}\psi_2)$. Notice however that 
$X(\psi_1,\psi_2)=-X(\psi_2,\psi_1)$.
\end{remark}

\bdeff
The vector field
\begin{equation}
\Xpao(\uu)=X(\phi_{j}(\uu),\phi_{j+1}(\uu)),
\end{equation}
defined in $\omega \setminus  \{\bar \uu\}$,  is called the \campopao\ associated with the conical intersection $\bar \uu$.
\edeff
The \campopao\ is smooth in its domain of definition. 
From~\eqref{obs2} and~\eqref{Xpaoprimo}, we have $\langle\phi_{j},\dot\phi_{j+1}\rangle = 0 $ along its integral curves.
Moreover, a simple computation leads to
\begin{align}\label{determinante}
 X(\psi_1,\psi_2)\cdot
\begin{pmatrix}
\langle \psi_2, H_1 \psi_2\rangle-\langle \psi_1, H_1 \psi_1\rangle\\
\langle \psi_2, H_2 \psi_2\rangle-\langle \psi_1, H_2 \psi_1\rangle\\
\langle \psi_2, H_3 \psi_2\rangle-\langle \psi_1, H_3 \psi_1\rangle
\end{pmatrix}
&\!=\! \frac{1}{2i} \det \mathcal{M}(\psi_1,\psi_2),
\end{align}
which implies the following result.
\bp\label{diff}
Let $\gamma(\cdot)\subset\omega\setminus\{\bar\uu\}$ be an integral curve of the \campopao. Then
\[\frac{d}{dt}\big[\lambda_{j+1}(\gamma(t))-\lambda_{j}(\gamma(t))\big]=\frac1{2i} F(\gamma(t)).\]
In particular, 
all the integral curves of the \campopao\ starting from a punctured neighborhood of the conical intersection reach it in  finite time (up to a time reversal).
\ep
 
 Without loss of generality, we can assume that $\bar \uu=0$.
Denote every $\uu\in \mathbb{R}^3$ by $\uu=\rho\vv$, where $\rho$ is its module and $\vv$ the versor. Let us call
$\mu(\rho,\vv)=\mm(\phi_j(\rho,\vv),\phi_{j+1}(\rho,\vv))$, where $\phi_j(\rho,\vv),\phi_{j+1}(\rho,\vv)$
denotes a choice of the eigenstates relative respectively to $\lambda_{j}(\rho\vv)$ and $\lambda_{j+1}(\rho \vv)$  (extended by continuity for $\rho=0$). We remark that $\mu(\rho,\vv)$ is defined up to a phase, therefore its module and the quantity
\begin{equation} \label{prodotto vet}
\mu(\rho,\vv)\times \mu^*(\rho,\vv)= 2i \ \mathfrak{Im}(\mu(\rho,\vv))\times \mathfrak{Re}(\mu(\rho,\vv))
\end{equation}
are well defined.
Finally,  set $\Xpaot(\rho,\vv)=\Xpao(\rho\vv)$. Notice that, for different values of $\vv$,  $\Xpaot(\rho,\vv)$ has a priori different limits as $\rho\to 0$.

The following estimates hold.
\begin{lemma} \label{stima rho}
Assume that $\bar \uu=0$ is a conical intersection. Then the  inequality
%at page \pageref{prodotto vet}. Then 
\begin{equation}
|\mu(\rho,\vv)\cdot \vv|\leq C\rho \label{er batman}
\end{equation}
 holds in a neighborhood of $0$ for some constant $C>0$ uniform with respect to $\vv$.
\end{lemma} 

\proof
Without loss of generality, we assume that 
the double eigenvalue is equal to 0. Then 
\begin{align*}
\rho\mu(\rho,\vv)\cdot \vv&=\rho\,\langle\phi_j(\rho,\vv),(v_1 H_1+  v_2 H_2+v_3H_3)\phi_{j+1}(\rho,\vv)\rangle\\
&=-{\langle\phi_j(\rho,\vv),H_0\phi_{j+1}(\rho,\vv)\rangle}=-\langle \phi_j(\rho,\vv)-\boldsymbol{P}_0\phi_j(\rho,\vv),H_0 \phi_{j+1}(\rho,\vv) \rangle.
\end{align*}

Assumption {\bf (H0)} implies that
\begin{align*}
\|H_0\phi_{j+1}(\rho,\vv) \| &= \| \lam_{j+1}(\rho\vv)\phi_{j+1}(\rho,\vv)   - \sum_{i=1}^3 \rho v_i H_i\phi_{j+1}(\rho,\vv) \|\nonumber\\
&\leq  | \lam_{j+1}(\rho,\vv)|  + \rho \sum_{i=1}^3  \|H_i\phi_{j+1}(\rho,\vv)\| \nonumber\\
%&\leq  | \lam_{j+1}(\rho,\vv)|  \left(1+\alpha(\rho\vv) \rho \right)+\beta(\rho\vv) \rho  \nonumber\\
&\leq  c\rho,   
\end{align*}
for some $c>0$, locally around the intersection.
By smoothness of the projection, 
we get  that 
\[\big|\langle \phi_j(\rho,\vv)-\boldsymbol{P}_0\phi_j(\rho,\vv),H_0 \phi_{j+1}(\rho,\vv)\rangle\big|\leq C'\rho^2,\] 
for a suitable $C'>0$, 
then we get the thesis.
\eproof

\medskip

We are now ready to prove the following result, which provides some information on the behavior of the trajectories of the \campopao.

\bp \label{pseudo introduzione}
With the notations introduced above and for $\rho$ small enough, there exist three constants $c_1,c_2,c_3>0$ such that $c_1 \leq |\dot{\rho}|\leq c_2$ and $|\dot{\vv}|\leq c_3$ along the trajectories  of the \campopao . %{\color{red} on $\omega\setminus\{\uu\}$}
\ep
\proof
Direct  computations lead to the equations
\begin{align*}
\dot{\rho}&=\Xpaot \cdot \vv,\\
\dot{\vv}&=\frac{1}{\rho} \left( \Xpaot - 
\left(\Xpaot \cdot \vv\right) \vv \right).
\end{align*}
The upper bound for $|\dot{\rho}|$ comes easily from the $H(\uu)$-boundedness of $H_i$ for every $i$.

From~\eqref{determinante} and~\eqref{prodotto vet}, we get that $|\Xpaot|=|\mathfrak{Im}(\mu(\rho,\vv))\times \mathfrak{Re}(\mu(\rho,\vv))|\geq c$ for some $c>0$ in a neighborhood of the singularity. 
An immediate consequence is that the three vectors $\mathfrak{Re}(\mu(\rho,\vv)),\mathfrak{Im}(\mu(\rho,\vv))$ and $\Xpaot(\rho,\vv)$
are three linearly independent vectors in $\mathbb{R}^3$.
Then we define the corresponding three unit vectors $\zz_1(\rho,\vv)=\mathfrak{Re}(\mu(\rho,\vv))/|\mathfrak{Re}(\mu(\rho,\vv))|$,
$\zz_2(\rho,\vv)=\mathfrak{Im}(\mu(\rho,\vv))/|\mathfrak{Im}(\mu(\rho,\vv))|$ and $\ww(\rho,\vv) = \Xpaot/|\Xpaot|$. 
We can then decompose $\vv$ in this basis as 
\[
\vv=a(\rho,\vv) \zz_1(\rho,\vv) + b(\rho,\vv) \zz_2(\rho,\vv)+\gamma(\rho,\vv) \ww(\rho,\vv).
\]
Notice that 
there exists some $\delta<1$ such that $|\zz_1(\rho,\vv) \cdot \zz_2(\rho,\vv)|\leq \delta$ for every $(\rho,\vv)$. Estimating the scalar product of $\vv$
with $\zz_1$ and  $\zz_2$, using~\eqref{er batman}, we obtain that
\[
|a(\rho,\vv)+b(\rho,\vv)\, \zz_1(\rho,\vv)\cdot \zz_2(\rho,\vv)|\leq C\rho \qquad |b(\rho,\vv)+a(\rho,\vv) \, \zz_1(\rho,\vv)\cdot \zz_2(\rho,\vv)|\leq C\rho,
\]
that easily leads to $|a(\rho,\vv)|\leq \bar{C}\rho$ and $|b(\rho,\vv)|\leq \bar{C}\rho$, for some uniform constant $\bar C$.

From~\eqref{er batman} and the fact that $|\vv|=1$ we get that 
\begin{align*}
|\ww(\rho,\vv)\cdot \vv| &\geq|\gamma(\rho,\vv)\ww(\rho,\vv)\cdot \vv|\\
&=|(\vv-(a(\rho,\vv) \zz_1(\rho,\vv) + b(\rho,\vv) \zz_2(\rho,\vv)))\cdot \vv|\\
&=|1- |a(\rho,\vv) \zz_1(\rho,\vv) + b(\rho,\vv) \zz_2(\rho,\vv)|^2| \geq 1-4 \bar C\rho^2.
\end{align*}
This implies that 
\[
|\ww(\rho,\vv)-(\ww(\rho,\vv)\cdot\vv) \vv|^2=|\vv-(\ww(\rho,\vv)\cdot\vv) \ww(\rho,\vv)|^2=1-(\vv\cdot \ww(\rho,\vv))^2 \leq 8 \bar{C}\rho^2.
\]
The thesis comes from the fact that 
\[
\left| \left( \Xpaot - 
\left(\Xpaot \cdot \vv\right) \vv \right) \right| = |\Xpaot | |\ww(\rho,\vv)-(\ww(\rho,\vv)\cdot\vv) \vv|
\qquad
\mbox{and} 
\qquad
|\Xpaot \cdot \vv |=|\Xpaot | |\ww(\rho,\vv)\cdot\vv|.
\]

\eproof

\medskip

The following proposition is a generalization of~\cite[Proposition~5.9]{adia-nostro} in the three dimensional case. 
The proof follows the same lines, thanks to Proposition~\ref{pseudo introduzione}, and is thus omitted.
\bp
\label{surjpao}
For every unit vector ${\bf v}$ in $\R^3$ there exists an integral curve $\gamma : [-\eta,0] \rightarrow \omega$ of $\Xpao$ with $\gamma(0)=0$, $\eta>0$, such that
\[
\lim_{t\rightarrow 0^-} \frac{\dot{\gamma}(t)}{\|\dot{\gamma}(t)\|} ={\bf v}.
\]
\ep

Thanks to Proposition~\ref{pseudo introduzione}, the integral curves of the \campopao\  are $\con^1$ up to the singularity included. 
In particular, they satisfy the hypotheses of
Proposition~\ref{proj smooth} with $k=1$, so that  the projections $P_j(\uu)$ and $P_{j+1}(\uu)$ on the eigenspaces relative to the 
intersecting eigenvalues are $\con^1$ along the integral curves of the \campopao\  outside the singularity, and can be continuously extended  
at the singularity. On the other hand,
$\boldsymbol{P}_{\uu}$  is $C^1$ along such curves, singularity included.

We remark moreover that, if $\gamma:[t_0,t_1] \to \mathbb{R}^3$ is an integral curve of the \campopao\ 
such that $\lambda_j(\gamma(t))\neq \lambda_{j+1}(\gamma(t))$ for $t\in[t_0,t_1)$, by definition of the \campopao, it holds
\begin{equation} \label{non mixing}
P_j(\gamma(t))\dot{P}_{j+1}(\gamma(t))=0 \qquad P_{j+1}(\gamma(t))\dot{P}_{j}(\gamma(t))=0 \qquad \qquad \forall \ t\in [t_0,t_1).
\end{equation}
We have the following.
\bp \label{C1}
Along every integral curve of the \campopao, there is a choice of the eigenstates relative to the intersecting eigenvalues which is $\mathcal{C}^1$  up to the singularity included. 
\ep

\proof
Let $I=[-T,0]$, and let $\gamma : I \to \mathbb{R}^3$ be an integral curve of $\Xpao$  such that $\gamma(0)=\bar \uu$ is a conical intersection between $\lambda_j$ and $\lambda_{j+1}$.
Outside the singularity, the eigenstates are well defined, up to a phase.  To fix the phase, we 
set
\[
\psi_j(t)=\frac{P_j(\gamma(t)) \widehat{\psi}}{\|P_j(\gamma(t)) \widehat{\psi}\|},
\]
where $\widehat{\psi}$ is an eigenstate of $H(\bar \uu)$ relative to $\lambda_j(\bar \uu)=\lambda_{j+1}(\bar \uu)$ such that 
$\lim_{t\to 0}P_j(\gamma(t)) \widehat{\psi}\neq 0$. Possibly reducing $T$, we can assume that $P_j(\gamma(t)) \widehat{\psi}\neq 0$ on the whole $I$, thus $\psi_j(t)$ is a normalized eigenstate of $H(\gamma(t))$ relative to $\lambda_j(\gamma(t))$. Since $\|P_j(\gamma(t)) \widehat{\psi}\|\neq0$, in order to prove that $\psi_j(t)$ is 
$\mathcal{C}^1$ it is enough to prove that $P_j(\gamma(t)) \widehat{\psi}$ is.

Since $P_j(\gamma(t))+P_{j+1}(\gamma(t))=\boldsymbol{P}_{\gamma(t)}$ for $t\in[-T,0)$, and by~\eqref{non mixing}, we get that
$P_j(\gamma(t))\dot{P}_{j}(\gamma(t))=P_j(\gamma(t))\dot{\boldsymbol{P}}_{\gamma(t)}$, and therefore
\[
\langle \psi_j(t),\dot{P}_j(\gamma(t))\widehat{\psi}\rangle=  \frac{\langle \widehat{\psi},P_j(\gamma(t))\dot{\boldsymbol{P}}_{\gamma(t)} \widehat{\psi}\rangle}{\|P_j(\gamma(t)) \widehat{\psi}\|}
,
\]
which has limit for $t\to 0^-$.

Let us notice that $\boldsymbol{P}_{\gamma(t)}P_j(\gamma(t))\widehat{\psi}=P_j(\gamma(t))\widehat{\psi}$. This, together with~\eqref{non mixing}, implies that 
\begin{align*}
\dot{P}_j(\gamma(t))\widehat{\psi}&=\dot{\boldsymbol{P}}_{\gamma(t)} P_j(\gamma(t))\widehat{\psi}+\boldsymbol{P}_{\gamma(t)}\dot{P}_j(\gamma(t))\widehat{\psi}\\
&=\dot{\boldsymbol{P}}_{\gamma(t)}P_j(\gamma(t))\widehat{\psi}+\langle \psi_j(t),\dot{P}_j(\gamma(t))\widehat{\psi}\rangle\psi_j(t),
\end{align*}
where the right hand side has limit for $t\to 0^-$.

We can repeat the same procedure to show that there is a choice for $\psi_{j+1}(t)$ such that $\dot{\psi}_{j+1}(t)$ has limit for $t\to 0^-$.
\eproof

\medskip

As an immediate consequence of the above proposition, we get the following corollary.
\bcc
In a neighborhood of a conical intersection, the integral curves of its associated non-mixing field are $\con^2$ up to the singularity included. 
\ecc
The regularity results proved earlier can be improved as shown below.

\bp
In a neighborhood of a conical intersection, the integral curves of its associated non-mixing field are $\con^{\infty}$ up to the singularity included.
In particular, we can choose $\con^{\infty}$ eigenstates $\psi_j,\ \psi_{j+1}$ along such a curve, up to the singularity included.
\ep

\proof
Let $\gamma : [-T,0] \to \mathbb{R}^3$ be an integral curve of $\Xpao$  such that $\gamma(0)=\bar \uu$ and $|\dot{\gamma}(t)|>0$ for every $t\in[-T,0]$ (this is true up to choosing 
$T$ sufficiently small), and define the eigenstates 
$\psi_j(\cdot),\psi_{j+1}(\cdot)$ as in Proposition~\ref{C1}.
We prove by induction on $n$ that, if the curve is $\con^n([-T,0])$, then also
$\psi_j(\cdot)$ and $\psi_{j+1}(\cdot)$ are $\con^n([-T,0])$. In particular, this last fact implies that 
the integral curves of the \campopao\ are $\con^{n+1}([-T,0])$.

In Proposition~\ref{C1} the claim is proved for $n=1$. Assume that it is true up to $n-1$, with $n>1$, that is
$\psi_j(\cdot)$ and $\psi_{j+1}(\cdot)$ are $\con^{n-1}([-T,0])$, and then the integral curves of the \campopao\ are $\con^{n}([-T,0])$.
In particular, we know that $\boldsymbol{P}_{\gamma(\cdot)} \in \con^{n}([-T,0])$, and 
$P_j(\gamma(\cdot)), \ P_{j+1}(\gamma(\cdot)) \in \con^{n-1}([-T,0])$.

Differentiating $n$ times the identity $\psi_j(t)=\boldsymbol{P}_{\gamma(t)}\psi_j(t)$, 
we obtain that for every $t\in [-T,0)$ $\psi_j^{(n)}(t)$ is a linear combination of the terms $\boldsymbol{P}^{(l)}_{\gamma(t)}
\psi_j^{(n-l)}(t)$, for $l=0,\ldots,n$. By inductive hypothesis, all the terms relative to $l\geq1$ are known to be continuous on $[-T,0]$. Then we have to prove the 
continuity of $\boldsymbol{P}_{\gamma(t)}
\psi_j^{(n)}(t)$ at $t=0$, which is equivalent to the continuity of $\boldsymbol{P}_{\gamma(t)} P_{j}^{(n)}(\gamma(t))\widehat{\psi}$.

We can develop $\boldsymbol{P}_{\gamma(t)}P_j^{(n)}(\gamma(t))$ as
\begin{align*}
\boldsymbol{P}P_j^{(n)}&=P_jP_j^{(n)}+P_{j+1}P_j^{(n)}\\
&=P_j\boldsymbol{P}^{(n)}-P_jP_{j+1}^{(n)}+P_{j+1}P_j^{(n)},
\end{align*}
where we omitted to precise that  all projections are evaluated along $\gamma$ for simplicity of notation. Notice that the term $P_j\boldsymbol{P}^{(n)}$ is continuous on the closed interval.
By~\eqref{non mixing}, it follows 
$\frac{d^{n-1}}{dt^{n-1}}\left(P_j(\gamma(t))\dot{P}_{j+1}(\gamma(t))\right)\equiv 0$, $t\in[-T,0)$, then we can write $P_j(\gamma(t))P_{j+1}^{(n)}(\gamma(t))$
as a linear combination of the terms $P_j^{(l)}(\gamma(t))P_{j+1}^{(n-l)}(\gamma(t))$, for $l=1,\ldots,n-1$. Since all these terms are continuous  on $[-T,0]$
by induction hypothesis, it follows that also  $P_j(\gamma(t))P_{j+1}^{(n)}(\gamma(t))$ is.

Analogous computations prove the same for $P_{j+1}(\gamma(t))P_j^{(n)}(\gamma(t))$, therefore
also $\boldsymbol{P}P_j^{(n)}$ can be defined continuously on $[-T,0]$.
In particular, $\boldsymbol{P}_{\gamma(t)}\psi_j^{(n)}(t)$ is continuous on $[-T,0]$, and we get the thesis.
The same argument proves the smoothness of $\psi_{j+1}(t)$.
\eproof

\medskip

\begin{figure}
\begin{center}
\includegraphics[width=0.45\textwidth]{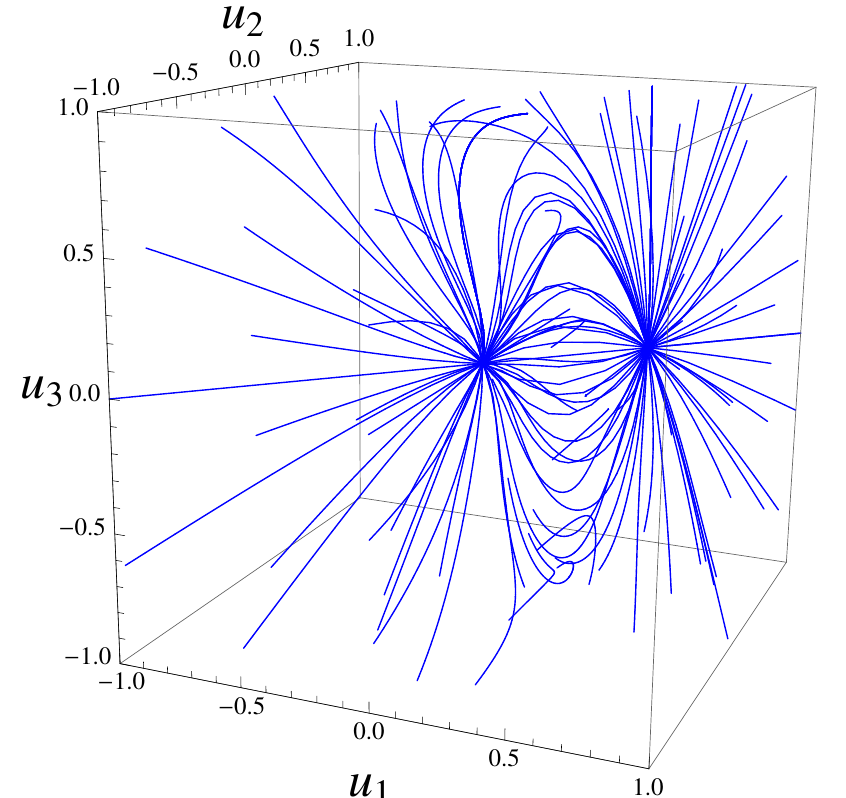}
\caption{Trajectories of the non-mixing field for the Hamiltonian corresponding to~\eqref{H-ricci}}\label{ricci}
\end{center}
\end{figure}

Figure~\ref{ricci} describes the behavior of the trajectories of the non-mixing field relative to the first two eigenvalues of the Hamiltonian~$H(\uu)\in i\mathfrak{u}(3)$
defined by \eqref{H-ricci}.
Consistently with the results shown above, the flow corresponding to the non-mixing vector field allows to identify two conical intersections, one of them being the origin,  among the two levels. 
In particular, the trajectories converge or 
diverge from them, locally.

\smallskip

To conclude this section, we present below a result providing some information on the structural stability of conical intersections based on the properties of the \campopao s.

\bt \label{struct}
Assume that $H(\uu)=H_0+u_1H_1+u_2H_2+u_3 H_3$ satisfies {\bf (H0)}-{\bf (H1)} and let $\bar\uu$ be a  conical intersection for $H(\uu)$ between the eigenvalues $\lambda_j$ and $\lambda_{j+1}$ belonging to the separated-discrete spectrum $\Sigma$. 
Then for every $\eps>0$ there exists $\delta>0$ such that, if $\hat{H}(\uu)=\hat H_0+u_1 \hat H_1+u_2\hat H_2+u_3\hat H_3$ satisfies  {\bf (H0)}-{\bf (H1)} and
\beq
\sum_{i=0}^3 \|\hat H_i-H_i\|_{H_0} \leq\delta,\label{varH}
\eeq
then the operator $\hat{H}(\uu)$ admits a conical intersection of eigenvalues at $\hat \uu$, with $|\bar \uu-\hat \uu|\leq \eps$. 
\et
\proof 
First of all, by equivalence of all norms $\|\cdot\|_{H(\uu)}$, without loss of generality we can assume that $\bar \uu=0$.
We notice that our assumptions guarantee that in a neighborhood of the conical intersection  the eigenvalues $\lambda_j$ and $\lambda_{j+1}$ are well separated from the rest of the spectrum.
Continuous dependence of the eigenvalues with respect to perturbations of the Hamiltonian (see Lemma~\ref{lemma-eigenva}) ensures that, if $\delta$ is small, then $\hat H(\cdot)$ admits two eigenvalues $\hat\lam_j,\hat\lam_{j+1}$ close 
to $\lam_j,\lam_{j+1}$.
Moreover $\{\hat\lam_j,\hat\lam_{j+1}\}$ is separated from the rest of the spectrum, locally around $\bar\uu$.

From the conicity of the intersection between $\lambda_j$ and $\lambda_{j+1},$
there exists $\eps>0$ small enough such that $|F(\uu)| \geq c$ for some $c>0$ on $B(\bar\uu,\eps)$ and moreover, by Proposition~\ref{pseudo introduzione}, 
the vector field $\Xpao$ (up to a global sign) points inside the ball $B(\bar\uu,\eps)$ at every point of its boundary.
If $\delta$ is small enough then $\hat\lam_j\neq\hat\lam_{j+1}$ on $\partial B(\bar\uu,\eps)$ and the gap between the two eigenvalues can be assumed to be of order $\eps$.
Therefore we can define the conicity matrix $\hat{\mathcal{M}}$ associated with $\hat{H}(\cdot)$ and the function $\hat{F}(\uu)=\det\hat{ 
\mathcal{M}}(\psi_j(\uu),\psi_{j+1}(\uu))$,
where $\{\psi_j(\uu),\psi_{j+1}(\uu)\}$ is an orthonormal basis for the sum of eigenspaces 
relative to $\{\hat\lam_j(\uu),\hat\lam_{j+1}(\uu)\}$. 
Since the conicity matrix varies continuously with respect to the control operators and the vectors on which it is evaluated,
we can take $\delta$ small enough such that $|\hat{F}(\uu)|\geq c/2$ on $B(\bar\uu,\eps)$.
This allows us to define, whenever $\hat\lam_j\neq\hat\lam_{j+1}$,  the non-mixing field $\hat{\mathcal{X}}_P$ associated with $\hat H(\cdot)$ and corresponding to the band  $\{\hat\lam_j,\hat\lam_{j+1}\}$; thanks to Proposition~\ref{diff}, up to a time reversal the time derivative of $\hat\lam_{j+1}-\hat\lam_{j}$ along the integral curves of  $\hat{\mathcal{X}}_P$ is smaller than $-c/4$.
By Corollary~\ref{high}  if $\delta$ is small enough, then $\hat{\mathcal{X}}_P$  points inside $B(\bar\uu,\eps)$ at every point of $\partial B(\bar\uu,\eps)$.  

Any trajectory $\hat\gamma(\cdot)$ of $\hat{\mathcal{X}}_P$ starting from $B(\bar\uu,\eps)$ 
 remains inside $B(\bar\uu,\eps)$ in its interval of definition and reaches in finite time a point $\hat \uu$ corresponding to a double eigenvalue  $\hat\lam_j(\hat\uu) = \hat\lam_{j+1}(\hat\uu)$. 
 The conclusion follows from Proposition~\ref{iff conical}.
\eproof

\subsection{A spread controllability result} \label{prova}

In this section we show how the curves tangent to the \campopao\ 
allow to improve the performances of the control algorithm presented in Section~\ref{control-trivial}. 
The result is stated here below.

\bt
\label{asc-paolo} 
Let $H(\uu)=H_0+   u_1 H_1+u_2 H_2+u_3 H_3$ satisfy hypotheses {\bf (H0)}-{\bf (H1)}.
Assume that there exist conical intersections $\bar \uu_j\in\omega,\ j=0,\ldots,k-1$, 
between the eigenvalues $\lam_j,\lam_{j+1}$, with $\lambda_l(\bar \uu_j)$ simple if $l\neq j,j+1$.
Then, 
for every $\uu^s$ and $\uu^f$ such that the eigenvalues $\lambda_{l},\ l=0,\ldots,k$ are non degenerate
at $\uu^s$ and $\uu^f$, for every
$\bar{\phi}\in \{\phi_0(\uu^s) ,\ldots,\phi_k(\uu^s) \}$, and $p \in [0,1]^{k+1}$ such that $\sum_{j=0}^k p_j^2=1$, 
there exist $C>0$ and a continuous control $\gamma(\cdot):[0,1]\to \R^m$ with $\gamma(0)=\uu^s$ and $\gamma(1)=\uu^f$, such that for every $\eps>0$
\begin{equation} \label{estimate}
\|\psi(1/\eps)-\sum_{j=0}^k p_j e^{i\vartheta_j} \phi_j(\uu^f)\|\leq C \varepsilon,
\end{equation}
where $\psi(\cdot)$ is the solution of (\ref{hamiltonianabold}) with $\psi(0)=\bar{\phi}$, $\uu(t)=\gamma(\eps t)$, and $\vartheta_0,\ldots,\vartheta_k \in \mathbb{R}$ are some phases depending on $\eps$ and $\gamma$.  
\et

\proof
The strategy is analogous to the one presented in Section~\ref{control-trivial} and consists in constructing a piecewise smooth path joining $\uu^s$ with $\uu^f$ that passes through all the conical intersection; in particular, we assume that the path  $\gamma : [0,1]\to \omega$ satisfies
$\gamma(0)=\uu^s$, $\gamma(1)=\uu^f$ and $\gamma(\tau_j)=\bar\uu_j, \ j=0,\ldots,k-1$, for some $0<\tau_0<\dots<\tau_{k-1}<1$. 
 The only difference concerns the choice of the paths in the neighborhoods of each conical intersection: indeed, in these regions our paths are chosen to be tangent 
 to the non-mixing field. At the intersection, the inner and outer directions are selected according to Proposition~\ref{intercono}, as explained in Section~\ref{control-trivial}, and the existence of corresponding trajectories  tangent to the non-mixing field is guaranteed by Proposition~\ref{surjpao}. 

Let us show that, given such a curve $\gamma(\cdot)$, the adiabatic approximations lead to the estimate~\eqref{estimate}.

Far from all the conical intersections, we 
approximate the evolution with that of the adiabatic Hamiltonian~\eqref{adiatutto}, which conserves the occupation probabilities relative to each energy level in the \mariospectrum .

On the other hand,  in a neighborhood of the conical intersection between $\lambda_j$ and $\lambda_{j+1}$, $j=0,\ldots,k-1$,  we approximate the dynamics with the 
ones associated with 
the adiabatic Hamiltonian~\eqref{adiaj}.

We show the details concerning the passage through $\bar \uu_0$, and the others can be treated analogously.
Since the path is tangent to the non-mixing field,  we can apply Theorem~\ref{effective} in order to study the evolution inside the 
space $\boldsymbol{P}_{\uu(t)}\H$.
For $\tau$ in a left neighborhood of $\tau_0$, we can then construct the effective Hamiltonian and its associated evolution operator $\Ueff$, which is diagonal, which implies that there exists a phase $\theta_0$ 
(depending on $\eps$) such that
\[
\|\psi(\tau_0/\eps)-e^{i\theta_0}\phi_0(\gamma(\tau_0^-))\| \leq C_0 \eps,
\]
where $\psi(\cdot)$ is the solution of equation~\eqref{hamiltonianabold} with $\psi(0)=\phi_0(\gamma(0))$ corresponding to the control $\uu(t)=\gamma(\eps t)$, defined on
$[0,1/\eps]$.
 
By Proposition~\ref{change of basis}, this implies that 
\[
\|\psi(\tau_0/\eps)-e^{i\theta_0}\left(\cos \bth(\vv)\ \phi_0(\gamma(\tau_0^+)) -e^{-i\beta(\vv)} \sin \bth(\vv)\ \phi_{1}(\gamma(\tau_0^+))\right)\| \leq C_0 \eps,
\]
with $\bth(\vv)$ and $\beta(\vv)$ satisfying equations~\eqref{alpha}-\eqref{beta}, and $\vv$ is the outer direction.

Since the effective Hamiltonian is diagonal also for $\tau$ belonging to a right neighborhood of $\tau_0$, we conclude that there exist two phases 
$\alpha_0$ and $\alpha_{1}$
(depending on $\tau$ and $\eps$) such that
\[
\|\psi(\tau/\eps)-e^{i\alpha_0}p_0 \ \phi_0(\gamma(\tau)) -e^{i\alpha_{1}} \sqrt{1-p_0^2} \ \phi_{1}(\gamma(\tau))\| \leq \hat{C}_0 \eps.
\]
Since analogous estimates hold on each passage through a conical intersection and outside the corresponding neighborhoods the theorem is proved.
\eproof

\section{Final remarks}
\label{fin rem}

The control strategy presented in the last section
(Theorem~\ref{asc-paolo}) highlights the role played by the integral curves of the non-mixing
field to obtain controllability results with an approximation of order $\eps$ on time intervals of order $1/\eps$ even in neighborhoods of conical intersections, while the classical theory would guarantee an error of order $\sqrt{\eps}$. 

As proved for the two-inputs case in~\cite{adia-nostro}, 
also for the three-inputs case it is possible to see that the 1-jet of the control path at the conical intersection determines the target probabilities, while the 2-jet is responsible of the
error introduced by the adiabatic approximation.
In other words,
$\con^2$ approximations (in the sense precised just below) of the integral curves of the \campopao\  also ensure an adiabatic
approximation of order $\eps$.
Indeed, let us consider an arc-length parametrized curve $\gamma^P(\tau)$, $\tau$ small, tangent to the non-mixing field,
that reaches a conical intersection between $\lambda_j$ and $\lambda_{j+1}$ at time $\tau=0$, and let
$\gamma$ be a $\con^3$ curve
such that $|\gamma(\tau)-\gamma^P(\tau)|\leq C \tau^3$, for $\tau$ small enough and some positive constant $C$. 
In particular, Proposition~\ref{proj smooth}
implies that there exists a constant $C'>0$ such that
 $|\phi_l(\gamma(\tau))-\phi_l(\gamma^P(\tau))|\leq C' \tau^2$ and
$|\dot{\phi}_l(\gamma(\tau))-\dot{\phi}_l(\gamma^P(\tau))|\leq C' \tau$,
$l=j,j+1$,
 where as usual $\phi_l(\uu)$ denotes the eigenstate relative to $\lambda_l(\uu)\in \Sigma(\uu)$,
evaluated at $\uu$. 
Let $t=\tau/\eps$ and consider the effective Hamiltonians evaluated along the two curves.
It is then easy to see by simple computations that there is a constant $C''>0$ such that the difference
between the two effective Hamiltonians is less or equal than  $C''\eps^2 t$. This term, integrated along a time
interval of order $1/\eps$, gives a difference of order $\eps$.

\smallskip

An interesting controllability problem alternative to the one introduced in Section~\ref{control-trivial} aims at sending (approximately) 
an initial state $\psi^s=\sum_{j=0}^k c_j\phi_j(\uu^s)$  
to a final state concentrated in a single energy level. This problem may appear 
completely equivalent to the previous one, but is actually more delicate, due to the presence of relative phases among the levels in the initial state. 
Indeed, a natural way to induce the desired transition would be to run backward in time one of the paths that produces any state with the same probability distribution as  $\psi^s$
starting 
from the concentrated state, constructed as in Section~\ref{control-trivial}  and Theorem~\ref{asc-paolo}. However, a simple computation shows that, 
at each passage through a conical intersection, the components corresponding to the intersecting eigenvalues recombine in a concentrated state only if  their relative phase coincides with the one induced by the unitary transformation of the limit basis (that is, the phase 
$\beta(\vv)$ of Proposition~\ref{change of basis}). 

There are several strategies that in principle could overcome this issue. For instance, when following the given path backward in time and before
reaching a conical intersection, it is always possible to stop for a certain time period at some point in the control space in order to control the relative
phase between the two intersecting levels.  Note that in the adiabatic evolution the effective stopping periods depend on the chosen speed $\eps$, 
while the geometric path in the space of controls does not depend on it. 
An alternative strategy consists in  exploiting the non-uniqueness  
property
underlined in Proposition~\ref{intercono}:  it is indeed possible to see that, for a path reaching a conical intersection and 
for any superposition of the two intersecting levels, there always exists a choice of outward direction allowing to concentrate the probability on a single level. Nevertheless, since the path is determined taking into account the dynamical phases, its construction depends on the time parameterization.

It is clear that the main drawback of the methods here above is that they rely on the computation of dynamical phases, which comes from the integration 
of the energy on intervals whose length is of order $1/\eps$ and which are very sensitive to changes in the speed $\eps$. 
This compromises the constructiveness of the algorithm. 

Let us finally mention that, to further improve the above controllability property, 
under assumptions {\bf (H0)}-{\bf (H1)} it is possible to modify the strategy, again with non-constructive arguments (for instance, by exploiting the rational 
independence of the gaps between the eigenvalues underlined in~\cite[Lemma 14]{ugauthier}) %see~\cite{ugauthier}),
in order to approximate at the final point not only any given choice of probability weights associated with $\Sigma$, but also any choice of the corresponding phases.

%On the other hand an important qualitative result which follows from 
Summing up, the constructions above combined with the control algorithm in Section~\ref{control-trivial} provide the following %(qualitative) 
approximate controllability result: 
\begin{quote}
under assumptions {\bf (H0)}-{\bf (H1)}, assuming that all energy levels in the separated discrete spectrum are connected through conical intersections, and
 for any given initial and target states $\psi^{s},\psi^{f}$ distributed in $\Sigma$ and $\eps>0$, 
 there exists a control input steering the system from $\psi^{s}$ to a final state whose distance from $\psi^{f}$ is less than $\eps$.
\end{quote}

\medskip
It is opinion of the authors that all the results here above still hold in a non-linear sufficiently smooth setting, that is for Hamiltonians of the form 
$H(\uu)$ whose derivatives with respect to the parameter $\uu$ are $H(\uu_0)$-small up to a suitable order and under hypothesis {\bf (H1)}. This case is interesting, since it covers relevant physical models, such as those described by Hamiltonians with controlled electromagnetic potentials. This topic will be the subject of further studies by the authors.

\appendix

\section{Regularity properties} \label{reg proj}

Let $\mathcal{H}$ be a complex separable Hilbert space;  all operators in the following are assumed to be operators on $\mathcal{H}$.
In this section  we derive some regularity results on the eigenvalues and the eigenstates of self-adjoint operators with respect to the norm defined in
\eqref{norm}. The results here below - partially already known in literature (see for instance \cite{katino}) - are proved by classical means.

In the following, $\rho(A)$ denotes the resolvent set of the operator $A$ and $\sigma(A)$ its spectrum. The resolvent of $A$  in $\zeta\in  \rho(A)$ is denoted
by $R(A,\zeta)=(A-\zeta \mathrm{id})^{-1}$; we recall that it  is a bounded linear operator that  maps $\mathcal{H}$ into $\mathcal{D}(A)$, and that,
given two self-adjoint operators $A_1,A_2$ with the same domain, their resolvents satisfy the Second Resolvent Identity 
\begin{equation} \label{rel res}
R(A_2,\zeta)- R(A_1,\zeta)=R(A_1,\zeta)(A_1-A_2)R(A_2,\zeta).
\end{equation}

First of all let us state the following technical lemma, which will be largely used in the following. Its proof easily comes from the definition of $\|\cdot\|_A$ 
and is thus omitted.
 \begin{lemma} 
Let $A,B$ be self-adjoint operators with $B$ $A$-bounded  and $\zeta\in\rho(A)$.
Then the following inequality holds:
\begin{equation} \label{stima risolvente}
\|BR(A,\zeta)\|\leq \Big(1+(|\zeta|+1)\, \|R(A,\zeta)\|\Big)\|B\|_A.
\end{equation}
\end{lemma}

The following result shows that the resolvent set for a self-adjoint operator $A$ enjoys some continuity properties with respect to small perturbation in the space $\mathcal{L}(\mathcal{D}(A),\mathcal{H})$. %, as the following lemma states. 
The proof follows from the definition of resolvent set $\rho(A)$ and properties of the resolvent $R(A,\zeta)$.

\begin{lemma} 
Let $A_1$ be a self-adjoint operator, and assume $[\zeta_1,\zeta_2]\subset \rho(A_1)$ for some real $\zeta_1\leq \zeta_2$. 
Then there exists a $\delta>0$ such that if $\|A_1-A_2\|_{A_1}\leq \delta$, then $A_1$ and $A_2$ have the same domain and  $[\zeta_1,\zeta_2]\subset \rho(A_2)$. 
Moreover, the inequality
\begin{equation} \label{stima risolvente delta}
\|R(A_2,\zeta)-R(A_1,\zeta)\|\leq \delta C %\frac{\delta C}{1-\delta C}
\end{equation}
holds on $[\zeta_1,\zeta_2]$ for some constant $C$ depending on $\zeta_1,\zeta_2$ and $A_1$.
\end{lemma}

\proof
Let $\zeta\in\rho(A_1)$. If the operator $\ \mathrm{id}+(A_2-A_1)R(A_1,\zeta)\ $ is invertible then the resolvent $R(A_2,\zeta)$ is well defined and bounded, and satisfies
\[R(A_2,\zeta)=R(A_1,\zeta)\Big(\mathrm{id}+(A_2-A_1)R(A_1,\zeta)\Big)^{-1}.\]
Thus the thesis follows once proved that, for every $\zeta\in[\zeta_1,\zeta_2]$, we have $\|(A_2-A_1)R(A_1,\zeta)\|<\delta C'$ for $\|A_1-A_2\|_{A_1}\leq \delta$ with $\delta$ small enough and for some $C'>0$.
This fact is a consequence of the uniform boundedness of $\|R(A_1,\zeta)\|$ on $[\zeta_1,\zeta_2]$ (see~\cite{katino}) and from~\eqref{stima risolvente}, and  
the thesis holds with $C=2\max_{\zeta\in[\zeta_1,\zeta_2]}\|R(A_1,\zeta)\|\Big(1+(|\zeta|+1)\, \|R(A_1,\zeta)\|\Big)$. 
\eproof

 \medskip

Let   $\lambda\in\sigma(A)$ be an eigenvalue of the self-adjoint operator $A$. For every positively-oriented closed path $\Gamma\subset\mathbb{C}$ 
encircling $\lambda$, and not encircling any other element in $\sigma(A)$, the projection $P$ onto the eigenspace relative to $\lambda$ is given by
\[
P=-(2\pi i)^{-1}\oint_{\Gamma} R(A,\zeta)\;d\zeta.
\]
\bp \label{cont proj}
Let $A_1$ be a self-adjoint operator, and let $\lambda$ be  a simple %isolated 
eigenvalue of $A_1$ such that $\sigma(A_1)\cap [\lambda-g,\lambda+g]=\{\lambda\}$ for some $g>0$. 
Then for every $\epsilon>0$ there exists a $\delta>0$ depending on $g$ and on $|\lambda|$ such  that if $\|A_1-A_2\|_{A_1}\leq \delta$, then  
\begin{enumerate}[i)]
\item  $\sigma(A_2)\cap [\lambda-g,\lambda+g]$ is made of only one point $\mu$, which is a simple eigenvalue of $A_2$;
\item Calling $P^{A_1}_{\lambda}$ the projection onto  the eigenspace of $A_1$ relative to $\lambda$ and $P^{A_2}_{\mu}$ the projection onto  the eigenspace of $A_2$ relative to $\mu$, it holds
\[
\|P^{A_1}_{\lambda}-P^{A_2}_{\mu}\| \leq \epsilon.
\]
\end{enumerate}
\ep

\proof
From preceding lemma, for every  $\delta>0$ small enough, if $\|A_1-A_2\|_{A_1}\leq \delta$, then 
$\sigma(A_2)\cap [\lambda-g,\lambda+g]$ is contained in the interval $(\lambda-g,\lambda+g)$.
Let $\Gamma$ be the circle in the complex plane of radius $g$ centered at $\lambda$, and consider the projection $P^{A_2}=-(2\pi i)^{-1}\oint_{\Gamma} R(A_2,\zeta)\;d\zeta$. From~\eqref{stima risolvente delta} we obtain that
\begin{align*}
\left\| P^{A_1}_{\lambda} - P^{A_2} \right\| &\leq (2\pi)^{-1}\oint_{\Gamma} \left\|R(A_1,\zeta)-R(A_2,\zeta) \right\| \;d\zeta\\
&\leq g C\delta
\end{align*}
that is, $\| P^{A_1}_{\lambda} - P^{A_2} \|< 1$ for $\delta$ small enough, which easily implies that 
$\dim \mathrm{Range}(P^{A_2})=\dim \mathrm{Range}(P^{A_1}_{\lambda})=1$, therefore $\sigma(A_2)\cap [\lambda-g,\lambda+g]$ contains exactly one spectral point
$\mu$,
which is a simple eigenvalue for $A_2$ (see~\cite[Theorem XII.6]{reed_simon}).
 
In particular, $P^{A_2}=P^{A_2}_{\mu}$ is the projection on the eigenspace relative to $\mu$, and satisfies {\it ii)} for $\delta$ small enough.
\eproof

\begin{mcc}
\label{high}
Under the hypothesis of Proposition~\ref{cont proj}, for every $\epsilon>0$ there exists a $\delta>0$ 
such  that if $\|A_1-A_2\|_{A_1}\leq \delta$, then  
\[
\|\phi^{A_1}_{\lambda}-\phi^{A_2}_{\mu}\| \leq \epsilon.
\]
where $\phi^{A_1}_{\lambda}$ and $\phi^{A_2}_{\mu}$ denote respectively the eigenstate of $A_1$ corresponding to $\lambda$ and
the eigenstate of $A_2$ corresponding to $\mu$ (normalized and with a particular choice for the global phases).
\end{mcc}

The next result provides an estimate concerning regularity properties  of the eigenvalues. 

\begin{lemma} 
\label{lemma-eigenva}
Let $A_1$ be a self-adjoint operator such that $\sigma(A_1)\cap I$ is discrete and without finite accumulation points for some open, possibly unbounded, interval $I$. If $\delta>0$ is small enough and $A_2$ is a self-adjoint operator satisfying $\|A_2-A_1\|_{A_1}\leq \delta$, then the eigenvalues of $A_2$ contained in $I$ are close to those of $A_1$,  in the following sense. 
Up to appropriately indexing on a subset   of $\mathbb{Z}$ the eigenvalues (counted with multiplicity) in $\sigma(A_j)\cap I$, for $j=1,2$,  and denoting them  with $\mu_i(A_j)$ we have $|\mu_i(A_1)-\mu_i(A_2)|\leq \epsilon (1+|\mu_i(A_1)|)$, where $\epsilon=e^{\frac{\delta}{1-\delta}} -1$.
\end{lemma}

\proof
Let $A_1$  satisfy the hypotheses of the lemma, and let $A_2$ be a self-adjoint operator with $\|A_1-A_2\|_{A_1}\leq \delta$, where  without loss of generality we assume 
that $\delta<1$;  define $A(t)=A_1+t(A_2-A_1)$, for $t\in[0,1]$. Let $\lambda_i(t)$ be the analytic branch of the eigenvalues of $A(t)$
emanating from $\lambda_i(0)=\mu_i(A_1)$, and denote by $\phi_i(t)$ a corresponding analytic eigenstate. 
By hypothesis
\begin{align*}
\|(A_2-A_1) \phi_i(t)\|  
&\leq \delta \|A_1 \phi_i(t)\| + \delta \\
&\leq \delta \|(A(t)-t(A_2-A_1)) \phi_i(t)\| + \delta \\
&\leq \delta |\lambda_i(t)|+\delta t\|(A_2-A_1))\phi_i(t)\| + \delta,
\end{align*}
which implies that $\|(A_2-A_1) \phi_i(t)\|  \leq \frac{\delta}{1-\delta } (|\lambda_i(t)|+1)$.

From 
\[
|\dot{\lambda}_i(t)|=|\langle \phi_i(t),(A_2-A_1)\phi_i(t)\rangle| \leq \|(A_2-A_1) \phi_i(t)\|  
\]
and Gronwall Lemma we easily get 
\[
|\lambda_i(0)-\lambda_i(1)| \leq \Big(e^{\frac{\delta}{1-\delta}} -1\Big)(|\lambda_i(0)|+1), 
\]
and then the thesis.
\eproof 

\medskip

%%%%%%%%%%%%%%%%%%%%%%

When we consider parameterized families of self-adjoint operators, we can prove some properties concerning the differentiability of spectral 
projections associated with separated portion of the spectrum. The statement here below deals with affine families, but  can be generalized to more general 
settings (see e.g.~\cite{katino,riesz-nagy} for similar arguments). 
\bp
\label{pdiff}
Let $K_0$ be a self-adjoint operator, $\mathcal{Y}$ be a Banach space with norm $\|\cdot\|_{\mathcal{Y}}$ and $K(\cdot)$ 
be a linear and continuous operator from $\mathcal{Y}$ to the space of $K_0$-bounded self-adjoint operators, endowed with the norm $\|\cdot\|_{K_0}$. Let moreover $q_0\in \mathcal{Y}$ and $I\subset \mathbb{R}$ be an interval whose boundary points belong to the resolvent set of $K_0+K(q_0)$. Then  the spectral projection $P_I(q)$ on $I$ associated with the self-adjoint operator $K_0+K(q)$
is well defined and (Fr\'echet) differentiable on a neighborhood of $q_0$.
\ep
\proof
By assumption we have that $\|K(q)-K(q_0)\|_{K_0} \leq C \|q-q_0\|_{\mathcal{Y}}$ for some $C>0$. Therefore %, thanks to~\eqref{stima risolvente}, 
for $q$ in a sufficiently small neighborhood of $q_0$, $\zeta$ belonging to the resolvent set of $K_0+K(q_0)$ and setting $R_{\zeta}(q)=R(K_0+K(q),\zeta)$, 
we can write, thanks to~\eqref{stima risolvente},
\[R_{\zeta}(q)=R_{\zeta}(q_0)\Big(\mathrm{id}+(K(q)-K(q_0))R_{\zeta}(q_0)\Big)^{-1}=R_{\zeta}(q_0)\sum_{k=0}^{\infty} \big((K(q_0)-K(q))R_{\zeta}(q_0)\big)^k.\]
Thus, from 
\[R_{\zeta}(q)-R_{\zeta}(q_0)-R_{\zeta}(q_0)(K(q_0)-K(q))R_{\zeta}(q_0)=R_{\zeta}(q_0) \sum_{k=2}^{\infty} \big((K(q_0)-K(q))R_{\zeta}(q_0)\big)^k\]
%{\color{Plum} and equation~\eqref{stima risolvente},}
we conclude that there exists a constant $\hat C>0$, continuously depending on $\zeta$, such that 
\[\|R_{\zeta}(q)-R_{\zeta}(q_0)-R_{\zeta}(q_0)(K(q_0)-K(q))R_{\zeta}(q_0)\|\leq \hat C \|q-q_0\|_{\mathcal{Y}}^2,\]
which guarantees that 
\[\lim_{\|q-q_0\|_{\mathcal{Y}}\to 0} \frac{1}{\|q-q_0\|_{\mathcal{Y}}} \left\| \oint_{\Gamma} \Big(R_{\zeta}(q)-R_{\zeta}(q_0)-R_{\zeta}(q_0)(K(q_0)-K(q))R_{\zeta}(q_0)\Big) d\zeta\right\|=0,\]
where $\Gamma$ is closed curve in $\mathbb{C}$ enclosing $I$ (and not containing any other element of $\sigma(K_0+  K(q_0))$).
\eproof

\medskip 
In the last part of this section, we focus on control-dependent Hamiltonians satisfying assumption {\bf (H0)} and, when explicitly said, 
assumption {\bf (H1)} too. 

First of all, it is easy to see that for any $\uu_1,\uu_2$ the norms $\|\cdot\|_{H(\uu_1)},\|\cdot\|_{H(\uu_2)}$ 
are equivalent, thanks to the $H_0$-smallness of the control Hamiltonians. 

Let us now focus on the eigenvalues of $H(\uu)$. From Lemma~\ref{lemma-eigenva}, 
and the equivalence of the norms $\|\cdot\|_{H(\uu)}$,
%under hypothesis {\bf (H0)} 
the eigenvalues $\lambda_i(\cdot)$ of $H(\cdot)$ are locally Lipschitz,
and the corresponding Lipschitz constants locally depend on the magnitude of $\lambda_i(\cdot)$. 
%and on $\uu$.

Let $\bar \uu$ be a conical intersection between the eigenvalues $\lambda_j$ and $\lambda_{j+1}$, that satisfy a gap condition, according to \textbf{(H1)}.
By the definition of conical intersection and the Lipschitz continuity of the eigenvalues we can conclude that there exist a suitably small neighborhood $U$ of $\bar \uu$  and two constants $C_1 > 0$ and $C_2>0$ such that
\begin{equation}\lambda_{j+1}(\uu) - \lambda_j(\uu) \geq C_1 |\uu - \bar \uu| \quad \forall \ \uu \in U  \label{concon}   \end{equation}
and
\begin{align} 
|\lambda_{i}(\uu) - \lambda_i( \uu')| &\leq C_2 |\uu -  \uu'| \quad \forall \ \uu,\uu' \in U,\ i=j,j+1. \label{Mj}
\end{align}

Moreover,  if we consider two eigenvalues $\lambda_j,\lambda_{j+1}$, possibly intersecting, and isolated from the rest of the spectrum, the projection $\boldsymbol{P}_{\uu}$ associated with these two levels is smooth with respect to $\uu$. The result holds also for any portion of the spectrum of $H(\uu)$, in presence of a gap (see~\cite{teufel}).

On the other hand, the projections $P_j$, $P_{j+1}$, associated respectively with $\lambda_j$ and $\lambda_{j+1}$, are smooth with respect to $\uu$ outside the singularity, while the presence of the conical intersection determines a lack of continuity at $\bar \uu$. Nevertheless, along regular curves passing through the singularity, it is possible to extend these projections, obtaining operators whose regularity depends on the one of the curve, as stated in the following result.

\bp \label{proj smooth}
Let $\gamma : I\to \mathbb{R}^3$, $I=[-R,0]$, be a $\mathcal{C}^k(I)$ curve such that $\gamma(0)=\bar \uu$ is a conical intersection between the eigenvalues $\lambda_j$ and $\lambda_{j+1}$
and 
$ \dot{\gamma}(t)\neq 0$ for every $t\in I$, and consider its k-jet at
the origin 
$\ell_k(t)=\gamma(0)+\sum_{j=1}^k \frac{1}{j!}t^j\frac{d^j}{dt^j}\gamma(t)|_{t=0}$.
Then $P_j(\gamma(\cdot))$ is $\mathcal{C}^{k}$ on $[-R,0)$, it is $\mathcal{C}^{k-1}$ at the singularity, and
\[
\lim_{t \to 0^-} \frac{d^{l}}{dt^{l}}P_j(\gamma(t))  =\lim_{t \to 0^-} \frac{d^{l}}{dt^{l}}P_j(\ell_k(t)),  \qquad l=0,\ldots,k-1,
\]
where the limit above holds in the operator norm.
%where we recall that the limit on the right-hand side of the above equation exists and is bounded, by virtue of the analiticity of $\ell_k$.
The same result holds for $P_{j+1}(\gamma(\cdot))$.
\ep

%%%%%%%%%%%%

\medskip
\noindent
%\textbf{Sketch of the proof. }
\proof
We first consider the case $k=1$. Without loss of generality we assume $|\dot{\gamma}(0)|=1$.
%Fix $\rho=\frac{1}{4}\min\{C_1,C_2\}$, 
Let $\rho=C_1/4$, where $C_1$ is as in~\eqref{concon}, and for every $t\in [-R,0)$ 
consider the circle $\Gamma_t\subset \mathbb{C}$ of radius $\rho t$ centered at  $\lambda_j(\gamma(t))$. 
For a set $A\subset \mathbb{C}$, we denote by $d(z,A)=\inf_{x\in A} |z-x|$ the distance between the point $z$ and the set $A$. 
There exists $0<T\leq R$ such that for every $t\in[-T,0)$ 
\[ |\lambda_{j+1}(\gamma(t))-\lambda_{j}(\gamma(t))| \geq \frac34 C_1 t = 3\rho t\]
so that $|\lambda_{j+1}(\gamma(t))-\zeta|\geq 2\rho t$ for every $\zeta\in \Gamma_t$. Thus  $d(\zeta,\sigma(H(\gamma(t))))=\rho t$ and, by~\eqref{Mj} and the definition 
of $\ell_1(\cdot)$, $d(\zeta,\sigma(H(\ell_1(t))))\geq \rho t/2$, up to reducing $T$.
Therefore from the classical identity holding for self-adjoint operators $\|(X-\zeta\mathrm{id})^{-1} \|=d(\zeta,\sigma(X))^{-1}$ (see e.g.~\cite{katino}), %pagina 277
for $\zeta\in \Gamma_t$ there hold 
\begin{equation}\label{stima risolv}
\|R(\gamma(t),\zeta)\| = \frac{1}{d(\zeta,\sigma(H(\gamma(t))))}=\frac{1}{\rho t}, \qquad 
\|R(\ell_1(t),\zeta)\| = \frac{1}{d(\zeta,\sigma(H(\ell_1(t))))}\leq \frac{2}{\rho t}.
\end{equation}

In order to get the thesis, we prove that 
\[
\lim_{t \to 0^-} \left\| \oint_{\Gamma_t} \left(R(\gamma(t),\zeta)- R(\ell_1(t),\zeta)\right)\; d\zeta \right\| =0.
\]
Estimate~\eqref{stima risolvente} gives 
\[
\|(H(\ell_1(t))-H(\gamma(t)))R(\gamma(t),\zeta)\| \leq C\left| \frac{\ell_1(t)-\gamma(t)}{t} \right| 
\]
for some $C>0$, which, together with~\eqref{rel res}-\eqref{stima risolv} and the definition of $\ell_1(\cdot)$, yields the thesis.

\smallskip
Let us now tackle the general case; the proof follows similar arguments. 
We define the circuit $\Gamma_{\tau}$ as above, and we notice that for every fixed $\tau\in (-T,0)$ there is a 
neighborhood $I_{\tau}$ of $\tau$ such that~\eqref{stima risolv} can be replaced by the similar estimate 
\begin{equation}\label{stima risolv 2}
\|R(\gamma(t),\zeta)\| \leq \frac{2}{\rho t}, \qquad 
\|R(\ell_k(t),\zeta)\| \leq \frac{2}{\rho t},
\end{equation}
holding for every $\zeta\in\Gamma_{\tau}$ and $t\in I_{\tau}$.
For every $t \in I_{\tau}$ we have that
\[
\frac{d^l}{dt^l}P_j(\gamma(t))=-(2\pi i)^{-1}\frac{d^l}{dt^l}\oint_{\Gamma_{\tau}} R(\gamma(t), \zeta)\;d\zeta
=-(2\pi i)^{-1}\oint_{\Gamma_{\tau}}\frac{d^l}{dt^l} R(\gamma(t), \zeta)\;d\zeta, \quad l\leq k-1 
\]
and thus, by applying~\eqref{rel res},
\[
\frac{d^l}{dt^l}P_j(\gamma(t))-\frac{d^l}{dt^l}P_j(\ell(t))=-(2\pi i)^{-1}\oint_{\Gamma_{\tau}}\frac{d^l}{dt^l} \Big(R(\ell_k(t),\zeta)(H(\ell_k(t))-H(\gamma(t)))R(\gamma(t),\zeta)\Big)\;d\zeta, \quad l\leq k-1. 
\]
The proof can then be easily completed by applying recursively the identity
\[\frac{d}{dt} R(f(t),\xi)=R(f(t),\xi)\Big(\frac{d}{dt}H(f(t))\Big)R(f(t),\xi),\]
the estimates~\eqref{stima risolv 2} and~\eqref{stima risolvente}, and by exploiting the regularity of $\gamma(\cdot)$ and the  definition of $\ell_k(\cdot)$.  %the regularity of the curve.
\eproof

%%%%%%%%%%%

\section{Genericity of conical intersections}\label{gen}

In this section we discuss the occurrence of conical intersection for certain classes of Hamiltonians. More precisely, we show that conical intersection
are not a pathological phenomena: indeed, in the finite dimensional case and for some specific families of infinite dimensional control-affine Hamiltonians, the fact
of being a conical intersection is a generic property of eigenvalue crossings, with respect to the controlled Hamiltonians.\footnote{
We recall that a property is said to hold generically in a Baire space $X$ if it is satisfied for all elements belonging to a residual subset of $X$, that 
is a set containing an intersection of (at most) countably many open and dense subsets of $X$.}

To prove this genericity property, in this paper we use transversality theorems that rely on the second countability of the family of Hamiltonian operators 
under consideration.
In the finite-dimensional case, all these hypotheses are fulfilled.
Concerning the general infinite-dimensional case, where we take the controlled Hamiltonians
$H_1,H_2,H_3$ as self-adjoint operators on $\mathcal{H}$, classical transversality theorems do not apply,
since the space of self-adjoint operators is
not second countable (even if we restrict our attention to bounded self-adjoint operators).
However physically relevant Hamiltonians often belong to particular families of operators which happen to be second-countable: we will focus on one of these 
cases.

In the following, we will consider the  class $\mathcal{F}$ of Hamiltonians of  the form $ \mathscr{K}(q)=K_0+K(q)$,  where $K_0$ is a fixed self-adjoint operator, $q$ belongs to a Banach space 
$\mathcal{Y}$  and $K(\cdot)$ is an injective linear and continuous operator from $\mathcal{Y}$ to the space of $K_0$-small self-adjoint operators, endowed with the norm $\|\cdot\|_{K_0}$.
We will also assume that the parameterized family $\mathcal{F}$ satisfies the following condition called Second Strong Arnold Hypothesis. %{\bf  (SAH2)}. 

\medskip
\noindent
\textbf{Second Strong Arnold Hypothesis (SAH2)} : {\it 
Assume that $\lambda$ is an eigenvalue of $\mathscr{K}(q_0)$ for some $q_0\in \mathcal{Y}$  of multiplicity greater or equal than two. Then there exist two 
orthonormal eigenstates
$\psi_1,\psi_2$ of $\mathscr{K}(q_0)$ pertaining to $\lambda$ such that the three linear functionals
\begin{align*}
f_{11} -f_{22} : q &\mapsto \langle \psi_1 , K(q) \psi_1\rangle -   \langle \psi_2 , K(q) \psi_2\rangle\\
f_{12} : q &\mapsto \langle \psi_1 , K(q) \psi_2\rangle\\
f_{21}  : q &\mapsto \langle \psi_2 , K(q) \psi_1\rangle
\end{align*}
are linearly independent.
Equivalently, the linear map $\Phi(\cdot)=(f_{11}(\cdot) -f_{22}(\cdot),\mathfrak{Re}(f_{12}(\cdot)),\mathfrak{Im}(f_{12}(\cdot)))$ is surjective from $\mathcal{Y}$ to $\mathbb{R}^3$.
}
\medskip

We call $\mathcal{D}$ the subset of $\mathcal{Y}$ such that the Hamiltonians in $\mathscr{K}(\mathcal{D})$ have double eigenvalues. 
For every interval $I$  and every open set $\mathcal{U}$ 
in $\mathcal{Y}$, we denote by $\mathcal{D}^{I,\mathcal{U}}$ the subset of elements in $\mathcal{U}$ such that the corresponding Hamiltonians have 
 an eigenvalue in $I$ of multiplicity two, isolated from the rest  of the spectrum.

In particular, under some additional regularity assumptions on the spectrum of the operators, {\bf  (SAH2)}  guarantees that $\mathcal{D}$  has codimension $3$ in $\mathcal{Y}$.
More precisely,  for a sufficiently small interval $I$ and a sufficiently small open set $\mathcal{U}$ the set $\mathcal{D}^{I,\mathcal{U}}$ 
is  a smooth manifold of codimension $3$ (see~\cite{teytel}).

The conicity of eigenvalue intersections correspond to a geometric property in the space of parameters, as the following result shows.
\begin{lemma} \label{tangente conica}
Let $H(\uu)=H_0+u_1 H_1+u_2 H_2+u_3 H_3$ belong to $\mathcal{F}$ for every $\uu \in \mathbb{R}^3$, 
%and satisfy {\bf (H0)-(H1)}, 
and assume that it has an isolated double eigenvalue $\lambda$ at $\uu=\bar \uu$; then there exists a unique $\bar{q} \in \mathcal{Y}$ such that $H(\bar \uu)=\mathscr{K}(\bar q)$
and  unique $q_1,q_2,q_3 \in \mathcal{Y}$ such that $H_i=K(q_i),\ i=1,2,3$.
Assume moreover that there is a neighborhood $\mathcal{U}$ of $\bar q$ in $\mathcal{Y}$ and an interval $I$ containing $\lambda$ such that 
$\mathcal{D}^{I,\mathcal{U}}$  is a submanifold of codimension three in $\mathcal{Y}$. 
Then $\bar \uu$ is a conical intersection 
for $H(\cdot)$ if and only if for every direction 
$\vv \in \mathbb{R}^3$ the vector $q_{\vv}=v_1q_1+v_2q_2+v_3q_3$ is not tangent
to $\mathcal{D}^{I,\mathcal{U}}$ at $\bar q$, that is, the affine space $\{\bar q+ q_{\vv} : \vv \in \mathbb{R}^3\}$ is transversal to $\mathcal{D}^{I,\mathcal{U}}$ at $\bar q$. 
%, that is, if and only if $H(\cdot)$ is transversal to $\mathcal{W}_d^{I,\mathcal{U}}$ at $\bar \uu$.
\end{lemma}

\proof
The existence of $\bar q$ and $q_i$ as in the thesis comes directly from linearity and injectivity of $K$.  

If there exists some $\vv \in \mathbb{R}^3$ such that $q_{\vv}$ is tangent
to $\mathcal{D}^{I,\mathcal{U}}$ at $\bar q$, then $\bar\uu$ cannot be a conical intersection, since in that case 
the distance between the eigenvalues intersecting at $\bar \uu$ is of  order $o(t)$ along the line $t\mapsto \bar\uu+t\vv$. 

Let us now prove the converse statement.
Denote by $\lambda_1(\uu)$ and $\lambda_2(\uu)$ the two eigenvalues
of $H(\uu)$ crossing at $\bar \uu$, with $\lambda_1(\bar\uu)=\lambda_2(\bar\uu)=\lambda$.

Under the assumptions of the lemma, we can deduce the following facts.
\begin{itemize}
 \item Possibly reducing $I$ (still containing $\lambda$ in its interior) and the neighborhood $\mathcal{U}$
of $\bar q$,  
$\mathscr{K}(q)$ contains exactly two eigenvalues in $I$, counted with their multiplicity, for every $q \in \mathcal{U}$.
\item Denoting  with $M(q)$ the sum of the eigenspaces of $\mathscr{K}(q)$ associated with the eigenvalues in $I$
and with $P_I(q)$ the orthogonal projection on $M(q)$, we have that, possibly reducing $\mathcal{U}$, 
\[\|P_I(q) - \bar P\| < 1 \qquad \forall q \in \mathcal{U},\]
where $\bar P=P_I(\bar q)$,
 and moreover $P_I(q)$ is a differentiable function of $q$ in $\mathcal{U}$, by Proposition~\ref{cont proj} and Proposition~\ref{pdiff}.
\item The map 
$S(q) = P_I(q) \big(\mathrm{id} + \bar P (P_I(q) - \bar P) \bar P \big)^{-1/2}\bar P$ 
is an isometric
transformation from $M(\bar q)$ onto $M(q)$ (see e.g.~\cite[Section~105]{riesz-nagy}), and it is differentiable with respect to its argument.
Therefore the map 
\[f(q) = S(q)^{-1}\mathscr{K}(q)S(q)\]
is a differentiable mapping from $\mathcal{U}$ to the space 
of self-adjoint operators on $M(\bar q)$, and the eigenvalues
of $f(q)$ are the same as the eigenvalues of $\mathscr{K}(q)$ in $I$.
\end{itemize}

It is easy to see that $\mathcal{D}^{I,\mathcal{U}}\subset f^{-1}(\{\mu \, \mathrm{id} : \mu\in \mathbb{R} \})$, 
where $\mathrm{id}$ denotes the identity 
on $M(\bar q)$. %matrix of dimension $2$.

Let us now assume that the  intersection between the eigenvalues $\lambda_1$ and $\lambda_2$ is not conical, 
that is there exists a unit vector $\mathbf{v}\in \mathbb{R}^3$ such that
\[
\lambda_2(\bar \uu+t \vv)-\lambda_1(\bar \uu+t \vv) =o(t),
\]
and we consider the curve $N(t)=f(\bar q + tq_{\vv})$ in the space of self-adjoint operators on $M(\bar q)$, that we write as two dimensional Hermitian 
matrices in a basis made of eigenstates relative to $\lambda_1$ and $\lambda_2$; 
it holds 
\[N(t)=
\begin{pmatrix}
a(t) & b(t)\\ b^*(t) & c(t) 
\end{pmatrix},\]
 for some complex functions $a(\cdot),b(\cdot),c(\cdot)$ satisfying $a(0)=c(0)=\lambda$ and $b(0)=0$.
Since the eigenvalues of $N(t)$ 
coincide with those of $H(\bar \uu+t\vv)$ contained in $I$, it holds $\sqrt{(a(t)-c(t))^2+4 |b(t)|^2}=o(t)$ and in particular,
by the analiticity of $\lambda_1(\bar \uu+t\vv)$ and $\lambda_2(\bar \uu+t\vv)$ with respect to $t$, it is easy to conclude that $\lim_{t \to 0^+}(\dot{a}(t)
-\dot{c}(t))=0$ and $\lim_{t \to 0^+}\dot{b}(t)=0$, that is, $N(t)$ is tangent to the space $\{\mu\, \mathrm{id} : \mu\in \mathbb{R} \}$ at the point $\lambda\, \mathrm{id}$. 

Since by definition the family $\mathcal{F}$ satisfies the condition {\bf (SAH2)}, the map $f$ is transversal to $\{\mu\, \mathrm{id} : \mu\in \mathbb{R} \}$. 
Then we can conclude that 
\[
(D_{\bar q}f)^{-1}(\{\mu\, \mathrm{id}: \mu\in \mathbb{R}\})=T_{\bar q}\mathcal{D}^{I,\mathcal{U}}
\]
(see e.g.~\cite[Corollary 17.2]{abraham-robbin}) and, in particular, $q_{\vv}$ is tangent to $\mathcal{D}^{I,\mathcal{U}}$ at $\bar q$.
\eproof

\subsection{Finite-dimensional case}
The class of Hamiltonians under consideration is here $\mathcal{F}=i\mathfrak{u}(n)$, that is  
 the set of Hermitian $n\times n$ matrices. Trivially, in this case $\mathcal{Y}$ coincides with $\mathcal{F}$. It is well
known (see e.g.~\cite{agrachev,von_neumann-wigner}) that
for sufficiently small $I$ and $\mathcal{U}$ the set $\mathcal{D}^{I,\mathcal{U}}$
is a \emph{smooth
manifold} of codimension 3 in $i\mathfrak{u}(n)$, and, moreover, that
the subset of Hermitian $n\times n$ matrices having at least one
eigenvalue
of multiplicity~3 has
codimension~8 in $i\mathfrak{u}(n)$. %Then, without loss of generality, we will focus on matrices having eigenvalues  of multiplicity at most two.
By second countability 
we can extract a countable family of pairs $\{(I_k,\mathcal{U}_k)\}_{k\in\mathbb{N}}$ that
satisfy the property above in such a way that $\mathcal{D}=\cup_k \mathcal{D}^{I_k,\mathcal{U}_k}$.

\begin{lemma}
\label{B52}
Fix $H_0 \in i\mathfrak{u}(n)$.
Borrowing notation from~\cite{abraham-robbin}, let us define the %following 
map $\rho: i\mathfrak{u}(n)^3 \to \mathcal{C}( \mathbb{R}^3\setminus\{0\}, i\mathfrak{u}(n))$
as
\begin{equation}
\rho(H_1,H_2,H_3)= H(\cdot) 
\end{equation}
with $H(\uu)=H_0+ u_1H_1+u_2H_2+u_3H_3$,
and let $\mathrm{ev}_{\rho} : (\mathbb{R}^3\setminus\{0\})\times i\mathfrak{u}(n)^3 \to i\mathfrak{u}(n)$ be defined as $\mathrm{ev}_{\rho} (\cdot, H_1,H_2,H_3)=\rho(H_1,H_2,H_3)$.
Then 
$\mathrm{ev}_{\rho}$ is transversal to $\mathcal{D}^{I_k,\mathcal{U}_k}$ for every $k\in \mathbb{N}$.
\end{lemma}

\proof
If $\mathrm{ev}_{\rho}(\uu,H_1,H_2,H_3)\notin \mathcal{D}^{I_k,\mathcal{U}_k}$, then the thesis trivially holds. 
Assume then that $H(\bar{\uu})=\mathrm{ev}_{\rho}(\bar{\uu},H_1,H_2,H_3)\in \mathcal{U}_k$ has a double eigenvalue $\lambda\in I_k$ for $\bar \uu\neq 0$.

The differential of $\mathrm{ev}_{\rho}$ at $(\bar{\uu},H_1,H_2,H_3)$ along the direction $(\delta\uu,\delta H_1, \delta H_2, \delta H_3)$ is given by
\[
D\mathrm{ev}_{{\rho}|_{(\bar{\uu},H_1,H_2,H_3)}} [\delta\uu,\delta H_1, \delta H_2, \delta H_3]= \sum_{i=1}^3 \left(\bar{u}_i \delta H_i+ \delta u_i H_i\right).
\] 
Let us consider the three directions
$v_l=(0,0,0,\bar{u}_1 \sigma_l,\bar{u}_2 \sigma_l,\bar{u}_3 \sigma_l)$, where the operators $\sigma_l$, $l=1,2,3$, are given by
\[
\sigma_1= \langle \varphi_2,\cdot \rangle\varphi_1 +  \langle \varphi_1,\cdot\rangle \varphi_2 \quad
\sigma_2=i \langle \varphi_1,\cdot\rangle\varphi_2 -i  \langle \varphi_2,\cdot\rangle \varphi_1
 \quad
\sigma_3= \langle \varphi_1,\cdot \rangle\varphi_1 -  \langle \varphi_2,\cdot\rangle\varphi_2, 
\]
and $\varphi_1$ and $\varphi_2$ define an orthonormal basis of the eigenspace of $H(\bar{\uu})$ relative to $\lambda$. 

Let us consider the eigenvalues of 
\[H(\bar{\uu})+\epsilon D\mathrm{ev}_{{\rho}|_{(\bar{\uu},H_1,H_2,H_3)}} [\alpha_1 v_1+\alpha_2 v_2+\alpha_3 v_3]=H(\bar{\uu})+\epsilon|\bar{\uu}|^2\sum_{i=1}^3 \alpha_i \sigma_i.\]
It is easy to check that the degenerate eigenvalues split and their difference is equal to
$2\epsilon |\bar{\uu}|^2 |\alpha|$. Therefore, $\mathrm{span}\{D\mathrm{ev}_{{\rho}|_{(\bar{\uu},H_1,H_2,H_3)}} [v_l] : l=1,2,3\}$  is a three dimensional space having trivial
intersection with $T_{H(\bar{\uu})}\mathcal{D}^{I_k,\mathcal{U}_k}$. 
Since the codimension of $\mathcal{D}^{I_k,\mathcal{U}_k}$ is 3, this means that $\mathrm{span}\{D\mathrm{ev}_{{\rho}|_{(\bar{\uu},H_1,H_2,H_3)}} [v_i] : i=1,2,3\}$ is transversal to 
$T_{H(\bar{\uu})}\mathcal{D}^{I_k,\mathcal{U}_k}$, 
and by definition of transversality of a map we get the proof.
\eproof

\bp \label{genericity finite}
Let $H_0\in i\mathfrak{u}(n)$. Then generically with respect to $(H_1,H_2,H_3)\in i\mathfrak{u}(n)^3$, all double eigenvalues of $H_0+\sum_{i=1}^3u_i H_i$ correspond to conical intersections.  
\ep

\proof
Define $H(\uu)=H_0+\sum_{i=1}^3 u_i H_i$. 
Thanks to Lemma~\ref{B52} we can apply the Transversal Density Theorem~\cite{abraham-robbin} with 
$\mathcal{W}=\mathcal{D}^{I_k,\mathcal{U}_k}$ for $k\in\mathbb{N}$, and $\rho(H_1,H_2,H_3)=H(\cdot)$, 
and obtain that the set of Hamiltonians $(H_1,H_2,H_3)$ such that $H(\cdot)$, as a function defined in $\mathbb{R}^3\setminus \{0\}$, is transversal to  
$\mathcal{D}^{I_k,\mathcal{U}_k}$ is residual in $i\mathfrak{u}(n)$.
In particular by applying Lemma~\ref{tangente conica} and taking the intersection of the previous residual sets over all $k$ we get that, generically, all double eigenvalues of $H(\uu)$ with $\uu\neq 0$ correspond to conical intersections.

Let us now consider the case where $H_0$ has a double eigenvalue $\lambda$, and let $\psi_1$ and $\psi_2$ define an orthonormal basis of the eigenspace relative 
to $\lambda$. 
Consider the real-valued multi-linear map
\[
F(H_1,H_2,H_3)=
\det\begin{pmatrix}
\mathfrak{Re}\left(\langle \psi_1, H_1 \psi_2 \rangle \right) & \mathfrak{Im}\left( \langle \psi_1, H_1 \psi_2 \rangle  \right)& \left( \langle \psi_2, H_1 \psi_2 \rangle-\langle \psi_1, H_1 \psi_1 \rangle \right)\\ 
\mathfrak{Re}\left(\langle \psi_1, H_2 \psi_2 \rangle \right) &\mathfrak{Im}\left( \langle \psi_1, H_2 \psi_2 \rangle \right) &\left(  \langle \psi_2, H_2 \psi_2 \rangle-\langle \psi_1, H_2 \psi_1 \rangle \right)\\
\mathfrak{Re}\left(\langle \psi_1, H_3 \psi_2 \rangle \right) &\mathfrak{Im}\left( \langle \psi_1, H_3 \psi_2 \rangle \right)& \left( \langle \psi_2, H_3 \psi_2 \rangle-\langle \psi_1, H_3 \psi_1 \rangle \right)
\end{pmatrix}.
\]
Recall that $\uu=0$ is a conical intersection for $H_0+\sum_{i=1}^3 u_i H_i$ if and only if $F(H_1,H_2,H_3)\neq 0$. Moreover, being $F$ continuous, we obtain that
$F^{-1}(\mathbb{R}\setminus\{0\})$ is an open subset of $i\mathfrak{u}(n)^3$.  This subset is non-empty because the map 
\[H\in i\mathfrak{u}(n)\mapsto
\Big(\mathfrak{Re}\left(\langle \psi_1, H \psi_2 \rangle \right),\  \mathfrak{Im}\left( \langle \psi_1, H \psi_2 \rangle  \right) ,\ \left( \langle \psi_2, H \psi_2 \rangle-\langle \psi_1, H \psi_1 \rangle \right)    \Big)   \in \mathbb{R}^3\]
is surjective.
 The density comes directly from multi-linearity.
\eproof

\subsection{Infinite dimension: the case of electromagnetic Hamiltonians}
\label{s-infinity}
Let us consider the class $\mathcal{F}$ of Hamiltonians of the form $H=-\Delta+V-i\left(\nabla \AV+\AV\nabla \right)$, where $\Delta$ 
denotes the Dirichlet Laplacian on a given Lipschitz bounded domain $\Omega\subset\mathbb{R}^3$, %\footnote{regularity assumptions?},
$V$ is a scalar continuous real-valued function on its closure $\bar\Omega$, %identified with 
that should be thought as a multiplication operator, and $\AV$ is a $\mathcal{C}^1$ vector-valued real function from $\bar\Omega$ to $\mathbb{R}^3$.
We focus on this family of Hamiltonians, since they happen to be largely used to model quantum systems driven by electromagnetic fields.

 The self-adjoint operator $H_{\AV}=-i\left(\nabla \AV+\AV\nabla \right)$  
acts on the elements of its domain 
as follows:
\[
-i\left(\nabla \AV+\AV\nabla \right)\psi=-i\, \AV \cdot \nabla \psi -i\,\mathrm{div}(\AV \psi). 
\]
Since, as it can be easily seen, the map $(V,\AV)\mapsto V+ H_{\AV}$ is linear and injective, $\mathcal{F}$ has a Banach manifold structure modeled on the space $\mathcal{Y}=\con(\bar\Omega,\mathbb{R}) \times 
\con^1(\bar\Omega,\mathbb{R}^3)$; moreover, since $\Omega$ is bounded, it turns out that both $\mathcal{Y}$ and $\mathcal{F}$ are separable, and thus 
second countable.

It is not difficult to show, by a direct integration by parts and applying the inequality $\langle \psi_1,\psi_2\rangle \leq \frac12
 (\eps \| \psi_1\|^2+\frac1{\eps}\|\psi_2\|^2)$, that self-adjoint operators of the form $H_{\AV}$ are $\Delta$-small. Therefore they can play 
the role of control Hamiltonians in our setting. Similarly it can be shown that each $H_{\AV}$ is form-bounded with respect to $-\Delta$ 
(as a quadratic form, see~\cite[Chapter~X]{reed_simon_2}) with a relative bound that can be chosen smaller than one. Thus~\cite[Theorem~XIII.68]{reed_simon} 
ensures that  the Hamiltonians of the form $-\Delta+V-i\left(\nabla \AV+\AV\nabla \right)$ have compact resolvent so that their spectrum 
is purely discrete with a finite number of eigenvalues in each compact subset of $\mathbb{R}$. Notice moreover that the topology induced by the norm $\|\cdot\|_{\Delta}$  on 
$\mathcal{F}$ is coarser than the topology  inherited from $\mathcal{Y}$.
 
For the family of Hamiltonians $\mathcal{F}$ defined above we essentially repeat the same argument as in the finite-dimensional case to show a  genericity property of conical intersections.
Before stating the main results of this section, some important remarks are in order. First of all we observe that the operator $H_{\AV}$ plays a crucial role 
for the existence of conical intersections for controlled Hamiltonian operators of the form $H_0+u_1H_1+u_2H_2+u_3H_3$ belonging to $\mathcal{F}$. 
Indeed if one considers controlled operators belonging to the class of Schr\"odinger operators  of the form $-\Delta+V$ with $V$ a real-valued function, 
then conical intersections are never present. This can be seen as a consequence of the fact that the %matrix elements 
terms $\langle \psi_j,H_i \psi_k\rangle$,
with $i=1,2,3$, computed with respect to an appropriately chosen orthonormal basis of eigenfunctions
$\{\psi_j\}_{j\in\N}$ of $-\Delta+V$, are real and thus the first two columns of each \mariomatrix\ are equal.
On the other hand, examples of conical intersections for controlled Hamiltonian belonging to the family $\mathcal{F}$ are not 
difficult to find, as shown below.

\medskip

\noindent{\bf Example.} Consider the Hamiltonian $H(\uu)=-\Delta+u_1V_1+u_2V_2+u_3H_{\AV}$, where 
\begin{gather*}
V_1(x)=x_2^2+x_3^2,\ V_2(x)=x_2x_3,\ A= (0,-x_3/2,x_2/2)^T,\\ 
x=(x_1,x_2,x_3)\in\Omega = (0,1)\times (0,\sqrt{3})\times (0,\sqrt{5})
\end{gather*}
with Dirichlet boundary conditions. We claim that $H(0)$ (representing the potential well in $\Omega$) admits conical intersections 
of eigenvalues. Indeed eigenvalues and eigenfunctions of $H(0)$ take the form
\[\lambda_{j_1,j_2,j_3}=\frac{\pi^2}{2} \Big(j_1^2+\frac{j_2^2}{3} +\frac{j_3^2}{5}\Big),\quad \psi_{j_1,j_2,j_3}(x) =\frac{2\sqrt{2}}{\sqrt[4]{15}} \sin ( j_1 \pi x_1)  \sin \Big( \frac{j_2 \pi x_2}{\sqrt{3}}\Big)\sin \Big( \frac{j_3 \pi x_3}{\sqrt{5}}\Big)\]
where $j_1,j_2,j_3$ are strictly positive integers.
Then one has that for instance $\lambda_{1,1,3}=\lambda_{1,2,2}$ corresponds to a double eigenvalue. A direct computation  shows that the associated \mariomatrix\ is nonsingular.

\medskip

Here we focus on the three-dimensional case since, unlike the other cases, it has a clear physical interest: in that case the vector $\AV$, called  vector potential, 
is related to the action of a magnetic field $\mathbf{B}$  on the physical system determined by the relation $\mathbf{B}=\nabla\times \AV$.
However, let us observe that the fact that the domain $\Omega$ is assumed to be a subset of $\R^3$ is not crucial for the mathematical 
formulation of the problem and the correctness of the following results  (the previous example, for instance, can be  directly recast 
in a two dimensional setting since the variable $x_1$ does not play any role there), unless the dimension is one. Indeed, in the latter case, 
it is easy to see that %, {\color{blue}in the one-dimensional case,}
for every continuously differentiable $\AV$, the energy levels of 
the Hamiltonian $-\partial_x^2 + V(x) + H_{\AV}$ coincide with those of the Hamiltonian  $-\partial_x^2 + V(x) - \AV^2$. 
 It can be easily shown that Hamiltonians of the latter form, on bounded intervals and with Dirichlet boundary condition,  do not admit degenerate eigenvalues. 
 Therefore the results below do not provide any information in the one dimensional case.
Note that the fact that for one dimensional systems a magnetic field can always be reabsorbed by an electric field is well known in physics.

Let us now proceed with the study of the genericity properties of conical intersections for the class $\mathcal{F}$ of controlled Hamiltonians under consideration. 

First of all, we notice that the class $\mathcal{F}$  fits the formulation given at the beginning of this appendix, with 
$\mathscr{K}(V,\AV)=-\Delta + V + H_{\AV}$.
We claim that the set $\mathcal{D}\subset \mathcal{Y}$ which parametrizes the operators admitting double eigenvalues  has codimension three in $\mathcal{Y}$.
The claim is proved once shown that all elements in $\mathcal{F}$ satisfy {\bf (SAH2)} (see~\cite{teytel}).

\begin{lemma} \label{SAH2}
Every element of $\mathcal{F}$ satisfies {\bf (SAH2)} for any  multiple eigenvalue. %In particular, $\mathcal{W}_d$ has codimension 3 in $\mathcal{F}$.
Moreover, the restriction $\Phi|_{\{(V,0) : V \in \con(\bar \Omega,\mathbb{R})\}}$ (defined in the statement of {\bf (SAH2)}) has rank at least two.
\end{lemma}

\proof
Let us consider $\bar H\in \mathcal{F}$, and assume that $\lambda$ is a multiple eigenvalue of $\bar H$.
By contradiction, assume that there exist three complex scalars $a,b,c$ and two eigenstates of $\bar H$ relative to the eigenvalue $\lambda$ such that the functional
\begin{equation} \label{funct}
 a (f_{11}(V,\AV) -f_{22}(V,\AV))+b f_{12}(V,\AV)+c f_{21}(V,\AV)
\end{equation}
is identically zero. 
%First of all, we notice that the nullness of the functional $\Phi$ does not depend on the 
Notice that this fact does not depend on the particular choice of the orthonormal eigenstates $\psi_1,\psi_2$.

Integrating by parts the terms of the kind $\langle \psi_i,H_{\AV}\psi_j\rangle$ %, and 
taking into account boundary conditions on the  eigenfunctions, %($\psi_i|_{\partial \Omega}=0$), 
we can write the functional above as
\[ 
\int_{\Omega} V \mathcal{B} - i \AV \cdot \mathcal{E},
\]
where
\begin{align*}
\mathcal{B}&=a |\psi_1|^2-a |\psi_2|^2 + b \psi_1^* \psi_2 + c \psi_1\psi_2^*, \\
 \mathcal{E}&=a  (\psi_1^* \nabla \psi_1- \psi_1 \nabla \psi_1^*- \psi_2^* \nabla \psi_2 + \psi_2 \nabla \psi_2^*)+ b 
 (\psi_1^* \nabla \psi_2 -\psi_2 \nabla \psi_1^*) + c (\psi_2^* \nabla \psi_1 -\psi_1\nabla \psi_2^*).
\end{align*}
By arbitrariness of $V$ and $\AV$, the expression~\eqref{funct} is identically zero only if $\mathcal{B}$ and $\mathcal{E}$ are identically zero on $\Omega$.

Let us start by assuming that there exist two orthonormal eigenfunctions $\psi_1,\psi_2$ such that $|\psi_1|\equiv |\psi_2|$ on $\Omega$ and three scalars $a,b,c$ 
such that the functional~\eqref{funct} is zero. Without loss of generality we can assume that $a\in \mathbb{R}$.
In particular, from $\mathcal{B}=0$ we obtain that either
$\psi_1$ and $\psi_2$ differ only by a constant phase, which contradicts their linear independence, or $b=c=0$. In the latter case, denoting 
$\psi_j=\psi e^{i\theta_j}$, $j=1,2$, it turns out that
$\mathcal{E}=2i a |\psi|^2 \nabla(\theta_1-\theta_2)$. Then $\nabla(\theta_1-\theta_2)=0$ wherever $|\psi|\neq 0$, and in particular
$\theta_1-\theta_2$ is constant on a open set. Up to a phase change of the eigenfunctions, 
$\psi_1-\psi_2$ is an eigenfunction which is null on a open set, which implies, by the unique
continuation property
(see e.g~\cite{kurata}), that
$\psi_1\equiv \psi_2$, which is a contradiction. 
We can then conclude that there are no orthonormal eigenfunctions with equal absolute value that make~\eqref{funct} identically zero.

Let us now assume the general case in which %that 
there exist  $\psi_1,\psi_2$  and  $a,b,c$ with $a\in\mathbb{R}$
such that the functional~\eqref{funct} is zero. %From
Again from $\mathcal{B}=0$, we obtain that 
\[
\mathfrak{Im}(\mathcal{B})=\mathfrak{Im}(b\psi_1^* \psi_2 + c \psi_1\psi_2^*)=(b-c^*)\psi_1^* \psi_2+(c-b^*)\psi_1 \psi_2^*\equiv 0 \quad \mbox{a.e. on } \Omega.
\]
If $b=c^*$, this condition is automatically satisfied and we write $\mathcal{B}$ as 
\[a |\psi_1|^2-a |\psi_2|^2 + b \psi_1^* \psi_2 + c \psi_1\psi_2^*=
 (\psi_1^*,\psi_2^*)
\begin{pmatrix}
a & b^*\\
b & -a
\end{pmatrix}
\begin{pmatrix}
 \psi_1\\ \psi_2
\end{pmatrix}.
\]
We can then diagonalize this quadratic form, ending up with two orthogonal eigenfunctions $\varphi_1,\varphi_2$ of $\bar H$, associated with $\lambda$, such that 
\[
\sqrt{a^2+|b|^2}(|\varphi_1|^2-|\varphi_2|^2)\equiv 0,
\]
 which is not possible thanks to the arguments above.

Let now $b\neq c^*$.
By unique continuation property, we can assume that $\psi_1^*\psi_2$ is not identically zero, that implies that
\[
\frac{b-c^*}{b^*-c}=\frac{\psi_1 \psi_2^*}{\psi_1^* \psi_2},
\]
which proves that the difference between the phases of $\psi_1$ and $\psi_2$ is constant on $\Omega$ and, in particular, it can be set to zero. This in particular 
leads to $\mathfrak{Im}(b)=-\mathfrak{Im}(c)$. Let us then set $b=\beta + i r$, $c=\gamma -i r$, for some $\beta,\gamma,r\in \mathbb{R}$.

Let us write $\psi_1=\phi_1 e^{i\zeta}$ and $\psi_2=\phi_2 e^{i\zeta}$, for some real-valued functions $\phi_1,\phi_2$ and $\zeta$.
Then by computations it follows from $\mathcal{E}=0$ that
\begin{align*}
&(\beta -\gamma)(\phi_1\nabla \phi_2 - \phi_2 \nabla \phi_1) + 2i\Big(r(\phi_1\nabla \phi_2 - \phi_2 \nabla \phi_1)+ a (\nabla \zeta) (\phi_1^2-\phi_2^2) 
+ (\gamma+\beta)(\nabla \zeta) \phi_1\phi_2\Big)=0.
\end{align*}
By direct computation one checks that $\phi_1\nabla \phi_2 - \phi_2 \nabla \phi_1$ is proportional to $\psi_1\nabla \psi_2 - \psi_2 \nabla \psi_1$,  and applying again the unique continuation property it turns out that the latter cannot be identically zero on an open set in $\Omega$. Therefore 
it must hold $\beta=\gamma$, that is $c=b^*$.
By contradiction, we see that {\bf (SAH2)} is verified.

The proof of the second statement follows the same arguments and is thus omitted.
\eproof

\medskip

We remark that, thanks to the fact that the spectrum of any operator in $\mathcal{F}$ is discrete with no finite accumulation points, and that the eigenvalues 
are continuous 
with respect to the pair $(V,\AV)$ (see Lemma~\ref{lemma-eigenva}), then for every $(\bar V,\bar \AV) \in \mathcal{Y}$ such that $\mathscr{K}(\bar V,\bar \AV)$ has a double eigenvalue $\lambda$ there exist a neighborhood $I$ of $\lambda$ and 
a neighborhood $\mathcal{U}$ of $(\bar V,\bar \AV)$ %in $\mathcal{Y}$ 
such that the subset $\mathcal{D}^{I,\mathcal{U}}$ is a smooth submanifold of codimension three in $\mathcal{Y}$. In particular, as in the finite dimensional case we can find a countable family $\{(I_k,\mathcal{U}_k)\}_k$ such that
$\mathcal{D}=\cup_k \mathcal{D}^{I_k,\mathcal{U}_k}$, with $\mathcal{D}^{I_k,\mathcal{U}_k}$ smooth  submanifold of $\mathcal{Y}$ of codimension three.

Let us now consider controlled Hamiltonians in $\mathcal{F}$ of the kind $H_0+u_1 V_1+u_2 V_2 +u_3 H_{\AV}$, where  $H_0\in \mathcal{F}$,
$V_1,V_2 \in \con(\bar \Omega,\mathbb{R})$ and $\AV \in \con^1(\bar \Omega,\mathbb{R}^3)$.

\begin{lemma} \label{F35}
Let $\bar H =- \Delta + \bar V + H_{\bar \AV}$, for some $\bar V \in \con(\bar \Omega,\mathbb{R})$ and $\bar \AV\in \con^1(\bar \Omega,\mathbb{R}^3)$.
Let $X=\{\uu \in \mathbb{R}^3 : u_1^2+u_2^2 \neq 0 \ \mathrm{and} \ u_3\neq 0\}$ and let us define the following map 
$\rho: (\con(\bar \Omega,\mathbb{R}))^2 \times \con^1(\bar \Omega,\mathbb{R}^3) \to \mathcal{C}(X , \mathcal{Y})$
as
\begin{equation}
\rho(V_1,V_2,\AV)(\uu)=(u_1 V_1+u_2 V_2, u_3 \AV). 
\end{equation}
Let $\mathrm{ev}_{\rho} : X\times (\con(\bar \Omega,\mathbb{R}))^2 \times  \con^1(\bar \Omega,\mathbb{R}^3) \to \mathcal{Y}$ 
be defined as $\mathrm{ev}_{\rho} (\cdot, V_1,V_2,\AV)=\rho(V_1,V_2,\AV)$.
Then $\mathrm{ev}_{\rho}$ is transversal to $\mathcal{D}^{I_k,\mathcal{U}_k}$ for every $k\in \mathbb{N}$.
\end{lemma}

\proof
Let us fix some notations: we set $K(V,\AV)=V+H_{\AV}$ and $\mathscr{K}(V,\AV)=\bar H + K(V,\AV)$.

If $\mathrm{ev}_{\rho} (\uu,V_1,V_2,\AV) \notin \mathcal{D}^{I_k,\mathcal{U}_k}$, then the thesis trivially holds. 

Assume then that $H(\bar{\uu})=\mathscr{K}(\mathrm{ev}_{\rho}(\bar{\uu},V_1,V_2,\AV))$ has a double eigenvalue $\lambda \in I_k$ 
for some $\bar \uu\in X$ and $(V_1,V_2,\AV)$ with $\mathrm{ev}_{\rho}(\bar{\uu},V_1,V_2,\AV) \in \mathcal{U}_k$, and that
$\psi_1,\psi_2$ are two orthonormal eigenstates of $H(\bar{\uu})$ pertaining to $\lambda$. Without loss of generality, we assume that $\bar u_1\neq 0$.
Lemma~\ref{SAH2} ensures that the rank of the map 
\[
(V,\AV) \mapsto \Big(\langle\psi_1, (V+H_{\AV})\psi_2 \rangle, \langle\psi_1, (V+H_{\AV})\psi_2 \rangle^*, 
\langle\psi_2, (V+H_{\AV})\psi_2 \rangle - \langle\psi_1, (V+H_{\AV})\psi_1 \rangle \Big)
\]
is three, and its restriction to the space $\{(V,0):V\in\con(\bar \Omega,\mathbb{R})\}$ has rank at least two. Then we can find three functions 
$\delta U,\delta W \in \con(\bar \Omega,\mathbb{R})$ and $\delta \AV \in \con^1(\bar\Omega,\mathbb{R}^3)$  such that
the conicity matrix 
\[
\mathcal{M}(\psi_1,\psi_2)=\begin{pmatrix}
\langle \psi_1,\delta U \psi_2 \rangle &  \langle \psi_1,\delta U \psi_2 \rangle^* &  \langle \psi_2,\delta U \psi_2 \rangle-\langle \psi_1,\delta U\psi_1 \rangle \\ 
\langle \psi_1,\delta W \psi_2 \rangle & \langle \psi_1,\delta W \psi_2 \rangle^* &  \langle \psi_2,\delta W \psi_2 \rangle-\langle \psi_1,\delta W \psi_1 \rangle \\
\langle \psi_1, H_{\delta \AV} \psi_2 \rangle & \langle \psi_1,  H_{\delta \AV} \psi_2 \rangle^* &  \langle \psi_2,  H_{\delta \AV} \psi_2 \rangle-\langle \psi_1,  H_{\delta \AV} \psi_1 \rangle
\end{pmatrix} 
\]
is non-singular.
In particular, this means that $\vv=0$ is a conical intersection for the Hamiltonian
\[
\widetilde{H}(\vv)=H(\bar \uu)+ v_1\delta U + v_2\delta W+ v_3 H_{\delta \AV}.\]
Moreover, we can write the Hamiltonian $\widetilde{H}(\vv)$ as 
\begin{align*}
\widetilde{H}(\vv) &=H(\bar \uu)+\sum_{j=1}^3 v_j K \big(  D\mathrm{ev}_{{\rho}|_{(\bar{\uu},V_1,V_2,\AV)}} [ w_j]\big) \\
&=\bar H + K(\bar{u}_1V_1+\bar{u}_2V_2,\bar{u}_3 \AV) + \sum_{j=1}^3 v_j K \big(  D\mathrm{ev}_{{\rho}|_{(\bar{\uu},V_1,V_2,\AV)}} [ w_j]\big)
\end{align*}
 where
$D\mathrm{ev}_{{\rho}|_{(\bar{\uu},V_1,V_2,\AV)}} [w]$ denotes
the differential of $\mathrm{ev}_{\rho}$ at $(\bar{\uu},V_1,V_2,\AV)$ evaluated on the variation 
$w \in \mathbb{R}^3 \times (\con(\Omega,\mathbb{R}))^2 \times  \con^1(\bar\Omega,\mathbb{R}^3)$, and the variations $w_j$ are 
\begin{gather*}
w_1=[0,0,0,\delta U/\bar{u}_1,0,0],\\ 
w_2=[0,0,0,\delta W/\bar{u}_1,0,0],\\ 
w_3=[0,0,0,0,0,\delta \AV/\bar{u}_3] .
\end{gather*}
Therefore, since the codimension of $\mathcal{D}^{I_k,\mathcal{U}_k}$ is three, and applying Lemma~\ref{tangente conica} with
$q_i=D\mathrm{ev}_{{\rho}|_{(\bar{\uu},V_1,V_2,\AV)}} [ w_i]$ and $\bar q=(\bar{u}_1V_1+\bar{u}_2V_2,\bar{u}_3 \AV) $, 
we get that 
\[\mathcal{Y}=\{q_{\vv} : \vv \in \mathbb{R}^3\} + T_{\bar{q}}\mathcal{D}^{I_k,\mathcal{U}_k}\subset\mathrm{Im}(D\mathrm{ev}_{{\rho}|_{(\bar{\uu},V_1,V_2,\AV)}}) + 
T_{\bar{q}}\mathcal{D}^{I_k,\mathcal{U}_k},\] 
and this concludes the proof.
\eproof

\bp \label{genericity infinite}
Let $H_0\in \mathcal{F}$ and $X$ as in the lemma above. Then generically with respect to $(V_1,V_2,\AV)\in (\con(\bar\Omega,\mathbb{R}))^2\times \con^1(\bar\Omega,\mathbb{R}^3)$, 
all the double eigenvalues of $H(\uu)=H_0+u_1 V_1 + u_2 V_2 +u_3 H_{\AV}$ with $\uu\in X$ correspond to conical intersections.  
\ep

\proof
Thanks to Lemma~\ref{F35} we can apply the Transversal Density Theorem 
and obtain that the set of triples 
$(V_1,V_2,\AV)$ such that $\rho(V_1,V_2,\AV)$, as a function defined on $X$,  is transversal to $\mathcal{D}^{I_k,\mathcal{U}_k}$ (at $\rho(V_1,V_2,\AV)$) is residual in $(\con(\bar\Omega,\mathbb{R}))^2\times \con^1(\bar\Omega,\mathbb{R}^3)$.
We can conclude as in Proposition~\ref{genericity finite}
that, generically, all double eigenvalues of $H(\uu)$ with $\uu\in X$ correspond to conical intersections.

\bigskip

\noindent{\it Acknowledgements.} 
This research has been supported  by the European Research Council, ERC
  StG 2009 ``GeCoMethods", contract number 239748, by the ANR ``GCM", program ``Blanc--CSD"
  project number NT09-504490, by the DIGITEO project ``CONGEO", and by the project CARTT-IUT Toulon.\\
The authors would like to thank Ugo Boscain for fruitful discussions.

\bibliographystyle{plain}
\bibliographystyle{abbrv}

\end{document}